\documentclass[12pt]{amsart}
\usepackage{amsmath}
\usepackage{amsfonts}
\usepackage{amssymb}
\usepackage[all,cmtip]{xy}           
\usepackage{mathrsfs}
\usepackage{bbm}
\usepackage{bbding}
\usepackage{txfonts}
\usepackage[shortlabels]{enumitem}
\usepackage{ifpdf}
\ifpdf
\usepackage[colorlinks,final,backref=page,hyperindex]{hyperref}
\else
\usepackage[colorlinks,final,backref=page,hyperindex,hypertex]{hyperref}
\fi
\usepackage{amscd}
\usepackage{tikz-cd}
\usepackage[active]{srcltx}

\makeatletter

\newtheorem{theorem}{Theorem}[section]
\newtheorem{prop}[theorem]{Proposition}
\newtheorem{lemma}[theorem]{Lemma}
\newtheorem{coro}[theorem]{Corollary}
\newtheorem{prop-def}{Proposition-Definition}[section]

\theoremstyle{definition}
\newtheorem{definition}[theorem]{Definition}

\newtheorem{remark}[theorem]{Remark}
\newtheorem{exam}[theorem]{Example}

\newcommand{\nc}{\newcommand}

\nc{\delete}[1]{{}}
\nc{\mmargin}[1]{}

\nc{\mlabel}[1]{\label{#1}}  
\nc{\mcite}[1]{\cite{#1}}  
\nc{\mref}[1]{\ref{#1}}  
\nc{\mbibitem}[1]{\bibitem{#1}} 

\delete{
	\nc{\mlabel}[1]{\label{#1}  
		{\hfill \hspace{1cm}{\bf{{\ }\hfill(#1)}}}}
	\nc{\mcite}[1]{\cite{#1}{{\bf{{\ }(#1)}}}}  
	\nc{\mref}[1]{\ref{#1}{{\bf{{\ }(#1)}}}}  
	\nc{\mbibitem}[1]{\bibitem[\bf #1]{#1}} 
}

 \font\cyrs=wncyr7

\nc{\vep}{\varepsilon}
\nc{\bin}[2]{ (_{\stackrel{\scs{#1}}{\scs{#2}}})}  
\nc{\binc}[2]{(\!\! \begin{array}{c} \scs{#1}\\
		\scs{#2} \end{array}\!\!)}  
\nc{\bincc}[2]{  ( {\scs{#1} \atop
		\vspace{-1cm}\scs{#2}} )}  
\nc{\oline}[1]{\overline{#1}}
\nc{\mapm}[1]{\lfloor\!|{#1}|\!\rfloor}
\nc{\bs}{\bar{S}}
\nc{\la}{\longrightarrow}
\nc{\ot}{\otimes}
\nc{\rar}{\rightarrow}
\nc{\lon }{\,\rightarrow\,}
\nc{\dar}{\downarrow}
\nc{\dap}[1]{\downarrow \rlap{$\scriptstyle{#1}$}}
\nc{\defeq}{\stackrel{\rm def}{=}}
\nc{\dis}[1]{\displaystyle{#1}}
\nc{\dotcup}{\ \displaystyle{\bigcup^\bullet}\ }
\nc{\hcm}{\ \hat{,}\ }
\nc{\hts}{\hat{\otimes}}
\nc{\hcirc}{\hat{\circ}}
\nc{\lleft}{[}
\nc{\lright}{]}
\nc{\curlyl}{\left \{ \begin{array}{c} {} \\ {} \end{array}
	\right .  \!\!\!\!\!\!\!}
\nc{\curlyr}{ \!\!\!\!\!\!\!
	\left . \begin{array}{c} {} \\ {} \end{array}
	\right \} }
\nc{\longmid}{\left | \begin{array}{c} {} \\ {} \end{array}
	\right . \!\!\!\!\!\!\!}
\nc{\ora}[1]{\stackrel{#1}{\rar}}
\nc{\ola}[1]{\stackrel{#1}{\la}}
\nc{\scs}[1]{\scriptstyle{#1}} \nc{\mrm}[1]{{\rm #1}}
\nc{\dirlim}{\displaystyle{\lim_{\longrightarrow}}\,}
\nc{\invlim}{\displaystyle{\lim_{\longleftarrow}}\,}
\nc{\dislim}[1]{\displaystyle{\lim_{#1}}} \nc{\colim}{\mrm{colim}}
\nc{\mvp}{\vspace{0.3cm}} \nc{\tk}{^{(k)}} \nc{\tp}{^\prime}
\nc{\ttp}{^{\prime\prime}} \nc{\svp}{\vspace{2cm}}
\nc{\vp}{\vspace{8cm}}
\nc{\modg}[1]{\!<\!\!{#1}\!\!>}
\nc{\intg}[1]{F_C(#1)}
\nc{\lmodg}{\!<\!\!}
\nc{\rmodg}{\!\!>\!}
\nc{\cpi}{\widehat{\Pi}}
\nc{\ssha}{{\mbox{\cyrs X}}} 
\nc{\tsha}{{\mbox{\cyrt X}}}
\nc{\shpr}{\diamond}    
\nc{\labs}{\mid\!}
\nc{\rabs}{\!\mid}

\nc{\C}{{\mathrm{C}}}
\nc{\ad}{\mrm{ad}}
\nc{\ann}{\mrm{ann}}
\nc{\Aut}{\mrm{Aut}}
\nc{\DA}{{\mathsf{DL}_\lambda}}
\nc{\Alg}{{\mathrm{Lie}}}
\nc{\DO}{{\mathsf{DO}_\lambda}}
\nc{\bim}{\mbox{-}\mathsf{Rep}}
\nc{\md}{\mbox{-}\mathsf{rep}}
\nc{\br}{\mrm{bre}}
\nc{\can}{\mrm{can}}
\nc{\Cont}{\mrm{Cont}}
\nc{\rchar}{\mrm{char}}
\nc{\cok}{\mrm{coker}}
\nc{\de}{\mrm{dep}}
\nc{\dtf}{{R-{\rm tf}}}
\nc{\dtor}{{R-{\rm tor}}}

\nc{\Div}{{\mrm Div}}
\nc{\Diff}{\mrm{DL}}
\nc{\Diffl}{\mathsf{DL}_\lambda}
\nc{\diffo}{{\mathsf{DO}_\lambda}}
\nc{\Dif}{{\mathfrak{Dif}^\lambda}}
\nc{\Difinfty}{{\mathfrak{Dif}^\lambda_\infty}}
\nc{\alg}{\mathsf{Lie}}
\nc{\End}{\mrm{End}}
\nc{\Ext}{\mrm{Ext}}
\nc{\Fil}{\mrm{Fil}}
\nc{\Fr}{\mrm{Fr}}
\nc{\Frob}{\mrm{Frob}}
\nc{\Gal}{\mrm{Gal}}
\nc{\GL}{\mrm{GL}}
\nc{\Hom}{\mrm{Hom}}
\nc{\Hoch}{\mrm{Hoch}}
\nc{\hsr}{\mrm{H}}
\nc{\hpol}{\mrm{HP}}
\nc{\id}{\mrm{id}}
\nc{\im}{\mrm{im}}
\nc{\Id}{\mrm{Id}}
\nc{\ID}{\mrm{ID}}
\nc{\Irr}{\mrm{Irr}}
\nc{\incl}{\mrm{incl}}
\nc{\length}{\mrm{length}}
\nc{\NLSW}{\mrm{NLSW}}
\nc{\Lie}{\mrm{Lie}}
\nc{\mchar}{\rm char}
\nc{\mpart}{\mrm{part}}
\nc{\ql}{{\QQ_\ell}}
\nc{\qp}{{\QQ_p}}
\nc{\rank}{\mrm{rank}}
\nc{\rcot}{\mrm{cot}}
\nc{\rdef}{\mrm{def}}
\nc{\rdiv}{{\rm div}}
\nc{\rtf}{{\rm tf}}
\nc{\rtor}{{\rm tor}}
\nc{\res}{\mrm{res}}
\nc{\Sh}{{\mathrm{Sh}}}
\nc{\SL}{\mrm{SL}}
\nc{\Spec}{\mrm{Spec}}
\nc{\sgn}{{\mathrm{sgn}}}
\nc{\tor}{\mrm{tor}}
\nc{\Tr}{\mrm{Tr}}
\nc{\tr}{\mrm{tr}}
\nc{\wt}{\mrm{wt}}
\nc{\op}{\mrm{op}}
\nc{\rmH}{ {\mathrm{H}}}

\nc{\bfk}{{\bf k}}
\nc{\bfone}{{\bf 1}}
\nc{\bfzero}{{\bf 0}}
\nc{\detail}{\marginpar{\bf More detail}
	\noindent{\bf Need more detail!}
	\svp}
\nc{\gap}{\marginpar{\bf Incomplete}\noindent{\bf Incomplete!!}
	\svp}
\nc{\FMod}{\mathbf{FMod}}
\nc{\Int}{\mathbf{Int}}
\nc{\Mon}{\mathbf{Mon}}
\nc{\remarks}{\noindent{\bf Remarks: }}
\nc{\Rep}{\mathbf{Rep}}
\nc{\Rings}{\mathbf{Rings}}
\nc{\Sets}{\mathbf{Sets}}
\nc{\ob}{\mathsf{Ob}}

\nc{\BA}{{\mathbb A}}   \nc{\CC}{{\mathbb C}}
\nc{\DD}{{\mathbb D}}   \nc{\EE}{{\mathbb E}}
\nc{\FF}{{\mathbb F}}   \nc{\GG}{{\mathbb G}}
\nc{\HH}{{\mathbb H}}   \nc{\LL}{{\mathbb L}}
\nc{\NN}{{\mathbb N}}   \nc{\PP}{{\mathbb P}}
\nc{\QQ}{{\mathbb Q}}   \nc{\RR}{{\mathbb R}}
\nc{\TT}{{\mathbb T}}   \nc{\VV}{{\mathbb V}}
\nc{\ZZ}{{\mathbb Z}}   \nc{\TP}{\widetilde{P}}
\nc{\m}{{\mathbbm m}}

\nc{\cala}{{\mathcal A}}    \nc{\calc}{{\mathcal C}}
\nc{\cald}{\mathcal{D}}     \nc{\cale}{{\mathcal E}}
\nc{\calf}{{\mathcal F}}    \nc{\calg}{{\mathcal G}}
\nc{\calh}{{\mathcal H}}    \nc{\cali}{{\mathcal I}}
\nc{\call}{{\mathcal L}}    \nc{\calm}{{\mathcal M}}
\nc{\caln}{{\mathcal N}}    \nc{\calo}{{\mathcal O}}
\nc{\calp}{{\mathcal P}}    \nc{\calr}{{\mathcal R}}
\nc{\cals}{{\mathcal S}}    \nc{\calt}{{\Omega}}
\nc{\calv}{{\mathcal V}}    \nc{\calw}{{\mathcal W}}
\nc{\calx}{{\mathcal X}}

\nc{\fraka}{{\mathfrak a}}
\nc{\frakb}{\mathfrak{b}}
\nc{\frakg}{{\frak g}}
\nc{\frakl}{{\frak l}}
\nc{\fraks}{{\frak s}}
\nc{\frakB}{{\frak B}}
\nc{\frakm}{{\frak m}}
\nc{\frakM}{{\frak M}}
\nc{\frakp}{{\frak p}}
\nc{\frakW}{{\frak W}}
\nc{\frakX}{{\frak X}}
\nc{\frakS}{{\frak S}}
\nc{\frakA}{{\frak A}}
\nc{\frakx}{{\frakx}}
\nc{\frakC}{{\frak{C}}}

\nc{\lir}[1]{\textcolor{red}{\underline{Li:}#1 }}

\begin{document}

\title [Twisted Rota-Baxter operators on 3-Leibniz algebras]{Twisted Rota-Baxter operators on 3-Leibniz algebras and NS-3-Leibniz algebras}

\author{Wen Teng}\footnote{Corresponding author}
\address{ School of Mathematics and Statistics, Guizhou University of Finance and Economics,  Guiyang  550025, P. R. of China}
\email{tengwen@mail.gufe.edu.cn}

\date{\today}

\begin{abstract}
The purpose of this paper is to introduce the  cohomology and deformations of twisted Rota-Baxter operators on 3-Leibniz algebras and NS-3-Leibniz algebras.
We construct an $L_\infty$-algebra whose Maurer-Cartan elements are twisted Rota-Baxter operators, and we define the cohomology of a twisted Rota-Baxter operator.
Then we consider formal and order $n$ deformations of twisted Rota-Baxter operators from cohomological points of view.
Finally,  we introduce and study NS-3-Leibniz algebras as the underlying structure of twisted Rota-Baxter operators on 3-Leibniz algebras.
\end{abstract}

\subjclass[2010]{
17A42, 17B38, 17B40, 17B56
}

\keywords{NS-3-Leibniz algebra;   twisted Rota-Baxter operator;  cohomology;  deformation;  $L_\infty$-algebra; Maurer-Cartan element}

\maketitle

\tableofcontents

\allowdisplaybreaks

\section{Introduction}

The class of $n$-Leibniz algebras, serving as a natural extension of Leibniz algebras to higher arities, was introduced in \cite{Casas}.
Furthermore, an $n$-Leibniz algebra, when combined with skew-symmetry, constitutes an $n$-Lie algebra \cite{Filippov}.
 In recent years, both domestic and international scholars have undertaken extensive research on the structure and properties of   $n$-Leibniz algebras.
 For instance,  in \cite{Albeverio},  the authors explored  the characteristics of Cartan subalgebras and normal
elements in $n$-Leibniz algebras. In  \cite{Casas16}, the authors employed Gr$\mathrm{\ddot{o}}$bner bases to formulate an algorithm for verifying the given multiplication table corresponding to an $n$-Leibniz algebra.
 In  \cite{Cherkis},  the authors developed three-dimensional $N=2$ superconformal Chern-Simons theories by imposing additional conditions on the 3-bracket.
 On this foundation, in \cite{Medeiros}, the authors investigated metric 3-Leibniz algebras and the hermitian 3-algebras introduced by Bagger-Lambert \cite{Bagger},
  which are intimately connected to $N=6$ and $N=5$ superconformal Chern-Simons theories.
Recently, in  \cite{XU}, the authors studied non-abelian extensions of 3-Leibniz algebras through Maurer-Cartan elements.
In  \cite{Hu24},  the authors introduced  3-pre-Leibniz algebras as the underlying structure of  relative
Rota-Baxter operators on 3-Leibniz algebras.
See  \cite{Azcarraga}  for more details about $n$-Leibniz algebras.

Rota-Baxter operators  were   introduced by Baxter \cite{Baxter60} in probability theory  and
were  further investigated  by Rota and Cartier  \cite{Cartier72,Rota69} in combinatorics.
 In order to understand the classical Yang-Baxter equation,  Kupershmidt \cite{{Kupershmidt99}} introduced the concept of a relative Rota-Baxter operator   on a Lie algebra.
After that, other operators related to the (relative) Rota-Baxter operators also appeared constantly. One of them is Reynolds operator, which is inspired by Reynolds's  work on  turbulence in fluid dynamics \cite{Reynolds95}.
In 1951,
Kamp$\mathrm{\acute{e}}$ de F$\mathrm{\acute{e}}$riet \cite{Kamp}   coined the term Reynolds operator, and then studied it as a mathematical subject.
For more interesting research on Reynolds operators, see \cite{CAlvarez,Chu,Rota,Zhang}.
Motivated by the twisted Poisson structure \cite{Severa}, Uchino \cite{Uchino} proposed twisted Rota-Baxter operators, which are the generalizations of  Reynolds operators and Rota-Baxter operators.
Recently, based on Uchino's work, Das  and his collaborators studied  the cohomology and deformation  theory of twisted Rota-Baxter operators on associative algebras, Lie algebras  and Leibniz algebras  in \cite{Das22,Das21,Guo24}.
 Furthermore,  the two papers \cite{Hou21,Chtioui}  independently introduced twisted Rota-Baxter operators on 3-Lie algebras, which are also called generalized Reynolds operators (twisted $\mathcal{O}$-operators) on 3-Lie algebras.


As a generalization of 3-Lie algebras, our primary objective in this paper is to investigate twisted Rota-Baxter operators on 3-Leibniz algebras and their associated structures.
To achieve this, we first  construct an $L_\infty$-algebra whose Maurer-Cartan elements are twisted Rota-Baxter operators on 3-Leibniz algebras. Moreover, we define a cohomology of a twisted Rota-Baxter operator
and apply it to control formal and order $n$  deformations of twisted Rota-Baxter operators.
Finally, we introduce a new ternary algebraic structure called NS-3-Leibniz algebras
that are related to twisted Rota-Baxter operators in the same way as 3-pre-Leibniz algebras are related to Relative Rota-Baxter operators.
A deeper exploration of NS-3-Leibniz algebras will be addressed in a forthcoming article.

The paper is organized as follows.  In Section \ref{sec: 3-Leibniz algebras}, we briefly recall basics about representations and
cohomology of 3-Leibniz algebras.
In Section \ref{sec: Twisted Rota-Baxter operators}, we introduce  twisted Rota-Baxter operators on 3-Leibniz algebras, provide some examples and characterization results.
Section \ref{sec: Maurer-Cartan} is devoted to constructing an $L_\infty$-algebra whose Maurer-Cartan elements are  twisted Rota-Baxter operators.
In Section \ref{sec: Cohomology}, we establish the cohomology theory of a  twisted Rota-Baxter operator  on a 3-Leibniz algebra.
 In Section \ref{sec: Deformations}, we study deformations of  twisted Rota-Baxter operators and show that they are controlled by the cohomology
theory established in Sections \ref{sec: Maurer-Cartan} and \ref{sec: Cohomology}. Finally, in Section \ref{sec: NS-3-Leibniz algebras} we introduce NS-3-Leibniz algebras which are derived naturally
from  twisted Rota-Baxter operators.

\section{ Basics on 3-Leibniz algebras }\label{sec: 3-Leibniz algebras}
\def\theequation{\arabic{section}.\arabic{equation}}
\setcounter{equation} {0}
In this section, we recall some basic notions about 3-Leibniz algebras and their representations  and cohomology theory.
We work over an algebraically closed field $\mathbb{K}$ of characteristic 0 and all the vector spaces
are over  $\mathbb{K}$.

\begin{definition} \cite{Casas}
A 3-Leibniz algebra is a vector space $\mathfrak{g}$ endowed with a  trilinear map $[\cdot, \cdot, \cdot]_\mathfrak{g}:\mathfrak{g}\otimes\mathfrak{g}\otimes\mathfrak{g}\rightarrow\mathfrak{g}$,
which satisfies:
\begin{align}
 &[a, b, [x, y, z]_\mathfrak{g}]_\mathfrak{g}=[[a, b, x]_\mathfrak{g},y,z]_\mathfrak{g}+ [x,  [a, b, y]_\mathfrak{g},z]_\mathfrak{g}+ [x,y,[a, b, z]_\mathfrak{g}]_\mathfrak{g} \label{2.1}
\end{align}
for all $ x, y, z, a, b\in \mathfrak{g}$.
\end{definition}

\begin{exam}
Any vector space with the trivial trilinear  bracket is a 3-Leibniz  algebra, called an abelian 3-Leibniz algebra.
\end{exam}

\begin{exam}
3-Lie algebras and Lie triple systems are 3-Leibniz algebras.
\end{exam}

\begin{prop} \label{prop: Leibniz-3Lie}
Let  $(\mathfrak{g},[\cdot,\cdot])$  be a Leibniz algebra. Then $(\mathfrak{g},[\cdot,\cdot,\cdot]_{\mathfrak{g}})$ is a 3-Leibniz algebra, where $[x,y,z]_\mathfrak{g}=[[x,y],z]$  for all
$x,y,z\in \mathfrak{g}$. Conversely,  if $(\mathfrak{g},[\cdot,\cdot,\cdot]_{\mathfrak{g}})$ is a 3-Leibniz algebra, then  $\mathfrak{g}\otimes\mathfrak{g}$ with the operation
$$[\mathcal{X},\mathfrak{Y}]=[x_1,x_2,y_1]_{\mathfrak{g}}\otimes y_2+ y_1\otimes[x_1,x_2,y_2]_{\mathfrak{g}}~\text{for all}~\mathcal{X}=x_1\otimes x_2,\mathfrak{Y}=y_1\otimes y_2\in \mathfrak{g}\otimes\mathfrak{g}$$
is a Leibniz algebra,
\end{prop}

\begin{definition} \cite{Casas} A representation of the 3-Leibniz algebra $(\mathfrak{g},[\cdot,\cdot,\cdot]_{\mathfrak{g}})$ is a
vector space $V$ equipped with three trilinear maps:
\begin{align*}
&\rho^l:\mathfrak{g}\otimes \mathfrak{g}\otimes V\rightarrow V,~~~\rho^m:\mathfrak{g}\otimes V\otimes \mathfrak{g}\rightarrow V,~~~\rho^r:V\otimes \mathfrak{g}\otimes \mathfrak{g}\rightarrow V,
\end{align*}
such that for any  $x,y,z,a,b\in \mathfrak{g}$ and $u\in V$,
\begin{align}
&\rho^l(a,b,\rho^l(x,y,u))=\rho^l([a,b,x]_{\mathfrak{g}}, y, u)+ \rho^l(x,[a,b,y]_{\mathfrak{g}},u)+\rho^l(x,y,\rho^l(a,b,u)),\label{2.2}\\
 &\rho^l(a,b,\rho^m(x,u,z))= \rho^m([a,b,x]_{\mathfrak{g}}, u, z)+\rho^m(x,\rho^l(a,b,u),z)+\rho^m(x,u,[a,b,z]_{\mathfrak{g}}),\label{2.3}\\
  &\rho^l(a,b,\rho^r(u,y,z))=\rho^r(\rho^l(a,b,u), y,z)+ \rho^r(u,[a,b,y]_{\mathfrak{g}},z)+\rho^r(u,y,[a,b,z]_{\mathfrak{g}}),\label{2.4}\\
   &\rho^m(a,u,[x,y,z]_{\mathfrak{g}})=\rho^r(\rho^m(a,u,x), y,z)+ \rho^m(x,\rho^m(a,u,y),z)+\rho^l(x,y,\rho^m(a,u,z)),\label{2.5}\\
    &\rho^r(u,b,[x,y,z]_{\mathfrak{g}})=\rho^r(\rho^r(u,b,x), y,z)+ \rho^m(x,\rho^r(u,b,y),z)+\rho^l(x,y,\rho^r(u,b,z)).\label{2.6}
\end{align}
\end{definition}

\begin{exam}
Given a 3-Leibniz algebra $(\mathfrak{g},[\cdot,\cdot,\cdot]_{\mathfrak{g}})$, there is a natural adjoint representation   on  itself.
The corresponding maps $\rho^l, \rho^m$ and $\rho^r$ are given by
$$\rho^l(x,y,z)=\rho^m(x,y,z)=\rho^r(x,y,z)=[x,y,z]_{\mathfrak{g}}  ~\text{for all}~  x,y,z\in \mathfrak{g}.$$
\end{exam}

\begin{prop} \label{prop:3-Leibniz} \cite{XU}
Let $(\mathfrak{g},[\cdot,\cdot,\cdot]_{\mathfrak{g}})$ be a 3-Leibniz algebra, $V$ be a vector space, and let
$\rho^l:\mathfrak{g}\otimes \mathfrak{g}\otimes V\rightarrow V,
\rho^m:\mathfrak{g}\otimes V\otimes \mathfrak{g}\rightarrow V $ and $\rho^r:V\otimes \mathfrak{g}\otimes \mathfrak{g}\rightarrow V$ be trilinear maps.
Then $(V;\rho^l,\rho^m,\rho^r)$  will be   a representation of  $\mathfrak{g}$ if and only if $(\mathfrak{g}\oplus V, [\cdot,\cdot,\cdot]_{\ltimes})$ is a  3-Leibniz algebra,
where $[\cdot,\cdot,\cdot]_{\ltimes}$ is defined as
\begin{align*}
[(x,u),(y,v),(z,w)]_{\ltimes}=&([x,y,z]_{\mathfrak{g}},\rho^l(x,y,w)+\rho^m(x,v,z)+\rho^r(u,y,z)),
\end{align*}
for all $(x,u),(y,v),(z,w)\in \mathfrak{g}\oplus V$.
\end{prop}

Let $(V;\rho^l,\rho^m,\rho^r)$ be a representation of a 3-Leibniz algebra  $(\mathfrak{g},[\cdot,\cdot,\cdot]_{\mathfrak{g}})$.
Denote by
$\mathcal{C}^n_{\mathrm{3Lei}}(\mathfrak{g},V)$ the set of n-cochains:
\begin{align*}
\mathcal{C}^n_{\mathrm{3Lei}}(\mathfrak{g},V)=\{\text{linear maps} ~f:\overbrace{\otimes^2\mathfrak{g}\otimes\cdots \otimes \otimes^2 \mathfrak{g}}^{n-1}\otimes \mathfrak{g}\rightarrow V,n\geq 1\}.
\end{align*}
 The coboundary map $\delta$
from $n$-cochains to
$(n+1)$-cochains, for $X_i=x_i\otimes y_i\in \otimes^2 \mathfrak{g}, 1\leq i\leq n$ and $z\in \mathfrak{g}$, as
\begin{align*}
&(\delta f)(X_1,X_2, \ldots, X_n,z)\\
=&\sum_{1\leq j<k\leq n}(-1)^jf(X_1,\ldots,\widehat{X_j},\ldots,X_{k-1},x_k \otimes[x_j,y_j,y_k]_{\mathfrak{g}}+[x_j,y_j,x_k]_{\mathfrak{g}}\otimes  y_k,\ldots,X_n,z)+\\
&\sum_{j=1}^n(-1)^jf(X_1,\ldots,\widehat{X_j},\ldots,X_{n},[x_j,y_j,z]_{\mathfrak{g}})+\sum_{j=1}^n(-1)^{j+1}\rho^l(X_j,f(X_1,\ldots,\widehat{X_j},\ldots,X_{n},z))+\\
&(-1)^{n+1}(\rho^m(x_n, f(X_1,\ldots,X_{n-1},y_n),z)+\rho^r(f(X_1,\ldots,X_{n-1},x_n), y_n, z)).
\end{align*}
It was proved in \cite{Daletskii} that $\delta^2=0$. Thus,
$(\oplus_{n=1}^{+\infty}\mathcal{C}^n_{\mathrm{3Lei}}(\mathfrak{g},V),\delta)$ is a cochain complex.  We denote the set of $n$-cocycles by
$\mathrm{Z}^n_{\mathrm{3Lei}}(\mathfrak{g},V)$, the set of $n$-coboundaries by $\mathrm{B}^n_{\mathrm{3Lei}}(\mathfrak{g},V)$ and the $n$-th cohomology group by $\mathrm{HH}^n_{\mathrm{3Lei}}(\mathfrak{g},V)=\frac{\mathrm{Z}^n_{\mathrm{3Lei}}(\mathfrak{g},V)}{ \mathrm{B}^n_{\mathrm{3Lei}}(\mathfrak{g},V)}$.

A trilinear map $\Phi\in\mathcal{C}^2_{\mathrm{3Lei}}(\mathfrak{g},V)$  is a 2-cocycle if $\delta \Phi=0$ which is equivalent to
\begin{align*}
0=(\delta \Phi)(x_1,y_1,x_2,y_2,z)=&-\Phi(x_2,[x_1,y_1,y_2]_{\mathfrak{g}},z)-\Phi([x_1,y_1,x_2]_{\mathfrak{g}},y_2,z)-\Phi(x_2,y_2,[x_1,y_1,z]_{\mathfrak{g}})+\\
&\Phi(x_1,y_1,[x_2,y_2,z]_{\mathfrak{g}})+\rho^l(x_1,y_1, \Phi(x_2,y_2, z))-\rho^l(x_2,y_2, \Phi(x_1,y_1, z))-\\
&\rho^m(x_2, \Phi(x_1,y_{1},y_2),z)-\rho^r(\Phi(x_1,y_{1},x_2), y_2, z).
\end{align*}
Under the above notations, the direct sum $\mathfrak{g}\oplus V$ carries a 3-Leibniz algebra structure
given by
\begin{align}
& [x+u,y+v,z+w]_\Phi=[x,y,z]_{\mathfrak{g}}+\rho^l(x,y,w)+\rho^m(x,v,z)+\rho^r(u,y,z)+\Phi(x,y,z), \label{2.7}
\end{align}
which is called the $\Phi$-twisted semi-direct product, denoted by $\mathfrak{g} \ltimes_\Phi V$.

At the end of this section, we give the relationship between the   representations of Leibniz algebra and 3-Leibniz algebra.
We recall Leibniz algebras and their representation.
A (left) Leibniz algebra is a vector space $\mathfrak{g}$ together with a bilinear operation
$[\cdot, \cdot]:\mathfrak{g}\times\mathfrak{g}\rightarrow\mathfrak{g}$ satisfying
\begin{align*}
 &[x, [y, z]]=[[x, y],z]+ [y,  [x,z]],~~\forall  x, y, z\in \mathfrak{g}.
\end{align*}

A representation of a Leibniz algebra $(\mathfrak{g},[\cdot,\cdot])$ consists of a triple $(V;\rho^L,\rho^R)$ in which
$V$ is a vector space and $\rho^L:\mathfrak{g}\otimes V\rightarrow V,$ $\rho^R:V\otimes\mathfrak{g}\rightarrow V$ are linear maps satisfying for any $x,y\in \mathfrak{g}$ and $u\in V$,
\begin{align*}
&\rho^L(x,\rho^L(y,u))=\rho^L([x,y],u)+\rho^L(y,\rho^L(x,u)),\\
 &\rho^L(x,\rho^R(u,y))=\rho^R(\rho^L(x,u),y)+\rho^R(u,[x,y]),\\
  &\rho^R(u,[x,y])=\rho^R(\rho^R(u,x),y)+\rho^L(x,\rho^R(u,y)).
\end{align*}
Notice that adding the last two equations, we have $\rho^R(\rho^L(x,u),y)+\rho^R(\rho^R(u,x),y)=0.$

\begin{theorem} \label{theorem:Nijenhuis}
Let   $(\mathfrak{g},[\cdot,\cdot])$ be a Leibniz algebra and $(V;\rho^L,\rho^R)$  be  a representation of  it. Then
$(V;\rho^l,\rho^m,\rho^r)$ is representation of 3-Leibniz algebra $(\mathfrak{g}, [\cdot,\cdot,\cdot]_\mathfrak{g})$  given in Proposition \ref{prop: Leibniz-3Lie},
where
\begin{align*}
&\rho^l(x,y,u)=\rho^L([x,y],u),~~\rho^m(x,u,y)=\rho^R(\rho^L(x,u),y)),~~\rho^r(u,x,y)=\rho^R(u,[x,y]),
\end{align*}
for all $x,y\in \mathfrak{g}$ and $u\in V$.
\end{theorem}
\begin{proof}
For any  $a,b,x,y,z\in \mathfrak{g}$ and $u\in V$, we have
  \begin{align*}
&\rho^l(a,b,\rho^l(x,y,u))=\rho^L([a,b],\rho^L([x,y],u))\\
&=\rho^L([[a,b],[x,y]],u)+\rho^L([x,y],\rho^L([a,b],u))\\
&=\rho^L([[a,b],x],y],u)+\rho^L([x,[[a,b],y]],u)+\rho^L([x,y],\rho^L([a,b],u))\\
&=\rho^l([a,b,x]_{\mathfrak{g}}, y, u)+ \rho^l(x,[a,b,y]_{\mathfrak{g}},u)+\rho^l(x,y,\rho^l(a,b,u)).
  \end{align*}
  Also,
  \begin{align*}
&\rho^l(a,b,\rho^m(x,u,z))=\rho^L([a,b],\rho^R(\rho^L(x,u),z))\\
&=\rho^R(\rho^L([a,b],\rho^L(x,u)),z)+\rho^R(\rho^L(x,u),[[a,b],z])\\
&=\rho^R(\rho^L([[a,b],x],u),z)+\rho^R(\rho^L(x,\rho^L([a,b],u)),z)+\rho^R(\rho^L(x,u),[[a,b],z])\\
&= \rho^m([a,b,x]_{\mathfrak{g}}, u, z)+\rho^m(x,\rho^l(a,b,u),z)+\rho^m(x,u,[a,b,z]_{\mathfrak{g}})
\end{align*}
and
\begin{align*}
&\rho^l(a,b,\rho^r(u,y,z))=\rho^L([a,b],\rho^R(u,[y,z]))\\
&=\rho^R(\rho^L([a,b],u),[y,z])+\rho^R(u,[[a,b],[y,z]])\\
&=\rho^R(\rho^L([a,b],u),[y,z])+\rho^R(u,[[[a,b],y],z])+\rho^R(u,[y,[[a,b],z]])\\
&=\rho^r(\rho^l(a,b,u), y,z)+ \rho^r(u,[a,b,y]_{\mathfrak{g}},z)+\rho^r(u,y,[a,b,z]_{\mathfrak{g}}).
\end{align*}
Similarly, we can prove the other Eqs. \eqref{2.5} and \eqref{2.6}.
  The proof is finished.
\end{proof}

\section{ Twisted Rota-Baxter operators on 3-Leibniz algebras }\label{sec: Twisted Rota-Baxter operators}
\def\theequation{\arabic{section}.\arabic{equation}}
\setcounter{equation} {0}
In this section we introduce twisted   Rota-Baxter operators on 3-Leibniz algebras.

\begin{definition}
A linear map $T:V\rightarrow \mathfrak{g}$ is called a $\Phi$-twisted Rota-Baxter operator on a 3-Leibniz algebra $(\mathfrak{g}, [\cdot,\cdot,\cdot]_\mathfrak{g})$ with respect to the representation  $(V;\rho^l,\rho^m,\rho^r)$ if $T$ satisfies
\begin{align}
& [Tu,Tv,Tw]_\mathfrak{g}=T(\rho^l(Tu,Tv,w)+\rho^m(Tu,v,Tw)+\rho^r(u,Tv,Tw)+\Phi(Tu,Tv,Tw)),  \label{3.1}
\end{align}
for all $u,v,w\in \mathfrak{g}$.
\end{definition}
\begin{exam}
Any relative Rota-Baxter operator  (in particular, Rota-Baxter operator of weight 0) on a 3-Leibniz algebra is a $\Phi$-twisted   Rota-Baxter operator with $\Phi=0$.
\end{exam}

Using the $\Phi$-twisted semi-direct product, one can characterize $\Phi$-twisted   Rota-Baxter operators by their graphs.

\begin{prop} \label{prop:graph}
A linear map $T:V\rightarrow \mathfrak{g}$  is a $\Phi$-twisted Rota-Baxter operator if
and only if its graph $Gr(T)=\{Tu+u~|~u\in V\}$ is a subalgebra of the $\Phi$-twisted semidirect
product $\mathfrak{g} \ltimes_\Phi V$.
\end{prop}

\begin{proof}
 Let $T:V\rightarrow  \mathfrak{g}$ be a linear map. Then, for all $u, v, w\in V$, we have
\begin{align*}
&[Tu+u, Tv+v, Tw+w]_{\Phi}\\
&=[Tu,Tv,Tw]_{\mathfrak{g}}+\rho^l(Tu,Tv,w)+\rho^m(Tu,v,Tw)+\rho^r(u,Tv,Tw)+\Phi(Tu,Tv,Tw),
\end{align*}
which implies that the graph  $Gr(T)$ is a  subalgebra of  $\mathfrak{g} \ltimes_\Phi V$
if and only if $T$ satisfies  Eq. \eqref{3.1}, which means that $T$ is  a  $\Phi$-twisted Rota-Baxter operator.
\end{proof}

Since $Gr(T)$  is isomorphic to $V$ as a vector space by the identification $(Tu,u)\backsimeq u$. A $\Phi$-twisted
Rota-Baxter operator $T$ induces a 3-Leibniz algebra structure on $V$ given by
\begin{align}
& [u,v,w]_T=\rho^l(Tu,Tv,w)+\rho^m(Tu,v,Tw)+\rho^r(u,Tv,Tw)+\Phi(Tu,Tv,Tw),  \label{3.2}
\end{align}
It is obvious that $T$ is a 3-Leibniz algebra morphism, that is $T[u,v,w]_T=[Tu,Tv,Tw]_\mathfrak{g}$.

\begin{definition}
Let $T:V\rightarrow \mathfrak{g}$  be a $\Phi$-twisted Rota-Baxter operator and $T':V'\rightarrow \mathfrak{g}'$  be a $\Phi'$-twisted Rota-Baxter operator.
A morphism of twisted  $\Phi$-twisted Rota-Baxter operators from $T$ to $T'$
consists of a pair $(f,g)$ of a 3-Leibniz algebra morphism
$f:\mathfrak{g}\rightarrow \mathfrak{g}'$
and a linear map $g:V\rightarrow V'$
satisfying
\begin{align}
&\left\{ \begin{array}{ll}
g(\rho^l(x,y,u))={\rho'}^{l}(f(x),f(y),g(u)),~g(\rho^m(x,u,y))={\rho'}^{m}(f(x),g(u),f(y)),\\
g(\rho^r(u,x,y))={\rho'}^r(g(u),f(x),f(y)),~g(\Phi(x,y,z))=\Phi'(f(x),f(y),f(z)),\\
f(Tu)=T'g(u),~~~~\forall ~~x,y,z\in \mathfrak{g},u\in V.
 \end{array}  \right.\label{3.3}
\end{align}
\end{definition}
\begin{exam}
Let $\mathfrak{g}$ be a 3-Leibniz algebra and  $(V;\rho^L,\rho^R)$  be  a representation of  it. Suppose $\varpi:\mathfrak{g}\rightarrow V$ is an invertible 1-cochain
in the cochain complex of  $\mathfrak{g}$  with coefficients in $V$. Then $T=\varpi^{-1}: V\rightarrow\mathfrak{g}$
is a $\Phi$-twisted Rota-Baxter-operator with $\Phi=-\delta \varpi$. To verify this, we observe that
\begin{align}
&{\Phi}(Tu,Tv,Tw)=-(\delta \varpi)(Tu,Tv,Tw) \label{3.4}\\
&=\varpi([Tu,Tv,Tw]_\mathfrak{g})-\rho^l(Tu,Tv,\varpi(Tw))-\rho^m(Tu,\varpi(Tv),Tw)-\rho^r(\varpi(Tu),Tv,Tw),\nonumber
\end{align}
By applying $T$ to the both sides of Eq. \eqref{3.4}, we get the Eq.  \eqref{3.1}.
\end{exam}

\begin{exam} \label{exam:Nijenhuis}
Let $N:\mathfrak{g}\rightarrow\mathfrak{g}$ be a Nijenhuis operator on a 3-Leibniz algebra $\mathfrak{g}$, i.e. $N$ satisfies
\begin{align}
[Nx,Ny,Nz]_\mathfrak{g}=&N\Big([x,Ny,Nz]_\mathfrak{g}+[Nx,y,Nz]_\mathfrak{g}+[Nx,Ny,z]_\mathfrak{g}\label{3.5}\\
&-N\big([x,y,Nz]_\mathfrak{g}+[x,Ny,z]_\mathfrak{g}+[Nx,y,z]_\mathfrak{g}\big)+N^2[x,y,z]_\mathfrak{g}\Big).\nonumber
\end{align}
In this case, $\mathfrak{g}$ carries a new  3-Leibniz algebra structure
\begin{align*}
[x,y,z]_N=&[x,Ny,Nz]_\mathfrak{g}+[Nx,y,Nz]_\mathfrak{g}+[Nx,Ny,z]_\mathfrak{g}\nonumber\\
&-N\big([x,y,Nz]_\mathfrak{g}+[x,Ny,z]_\mathfrak{g}+[Nx,y,z]_\mathfrak{g}\big)+N^2[x,y,z]_\mathfrak{g},
\end{align*}
We denote this 3-Leibniz algebra  by $\mathfrak{g}_N$. Moreover, the 3-Leibniz algebra $\mathfrak{g}_N$ has a representation on $\mathfrak{g}$ given by
$$\rho^l(x,y,z)=[Nx,Ny,z]_\mathfrak{g},\rho^m(x,y,z)=[Nx,y,Nz]_\mathfrak{g},\rho^r(x,y,z)=[x,Ny,Nz]_\mathfrak{g}.$$
With this representation, the map $\Phi:\otimes^3 \mathfrak{g}\rightarrow\mathfrak{g}$ defined by
$$\Phi(x,y,z)=-N\big([x,y,Nz]_\mathfrak{g}+[x,Ny,z]_\mathfrak{g}+[Nx,y,z]_\mathfrak{g}\big)+N^2[x,y,z]_\mathfrak{g}$$
is a 2-cocycle in the  cohomology of $\mathfrak{g}_N$ with coefficients in $\mathfrak{g}$. Then it is easy to
observe that the identity map $\mathrm{Id}:\mathfrak{g}\rightarrow\mathfrak{g}_N$ is a $\Phi$-twisted Rota-Baxter-operator.
\end{exam}

\begin{exam} \label{exam:Reynolds}
Let  $(\mathfrak{g}, [\cdot,\cdot,\cdot]_\mathfrak{g})$  be a 3-Leibniz algebra. Set $\Phi=-[\cdot,\cdot,\cdot]_\mathfrak{g}$;
 then, a linear transformation $T:\mathfrak{g}\rightarrow \mathfrak{g}$ defined by Eq. \eqref{3.1} is called a Reynolds operator;
more specifically, $T$ satisfies that
\begin{align*}
& [Tx,Ty,Tz]_\mathfrak{g}=T\big( [Tx,Ty,z]_\mathfrak{g}+ [Tx,y,Tz]_\mathfrak{g}+ [T,Ty,Tz]_\mathfrak{g}- [Tx,Ty,Tz]_\mathfrak{g}\big),
\end{align*}
for all $x,y,z\in \mathfrak{g}$.
\end{exam}

Given a $\Phi$-twisted Rota-Baxter-operator $T$ and a 1-cochain $\varpi$, we construct a $(\Phi+\delta\varpi)$-twisted Rota-Baxter-operator  under
certain condition. First we observe the following

\begin{prop} \label{prop:isomorphism of 3-Leibniz algebras}
Let $\mathfrak{g}$ be a 3-Leibniz algebra and  $(V;\rho^L,\rho^R)$  be  a representation of  it.
 For any 2-cocycle $\Phi\in\mathcal{C}^2_{\mathrm{3Lei}}(\mathfrak{g},V)$ and
1-cochain $\varpi\in\mathcal{C}^1_{\mathrm{3Lei}}(\mathfrak{g},V)$, we have an isomorphism of 3-Leibniz algebras
$$\mathfrak{g} \ltimes_\Phi V\cong\mathfrak{g} \ltimes_{\Phi+\delta\varpi} V.$$
\end{prop}

\begin{proof}
Define $\varphi_\varpi:\mathfrak{g} \ltimes_\Phi V\rightarrow\mathfrak{g} \ltimes_{\Phi+\delta\varpi} V$ by $\varphi_\varpi(x,u)=x+u-\varpi(x)$, for $x+u\in \mathfrak{g} \ltimes_\Phi V$. Then
 $x+u,y+v,z+w\in \mathfrak{g} \ltimes_\Phi V$, we have
\begin{align*}
&\varphi_\varpi[x+u,y+v,z+w]_\Phi\\
&= [x,y,z]_{\mathfrak{g}}+\rho^l(x,y,w)+\rho^m(x,v,z)+\rho^r(u,y,z)+\Phi(x,y,z)-\varpi([x,y,z]_{\mathfrak{g}})\\
&=[x,y,z]_{\mathfrak{g}}+\rho^l(x,y,w)+\rho^m(x,v,z)+\rho^r(u,y,z)+\Phi(x,y,z)+(\delta\varpi)(x,y,z)-\rho^l(x,y,\varpi(z))\\
&-\rho^m(x,\varpi(y),z)-\rho^r(\varpi(x),y,z)\\
&=[x+u-\varpi(x),y+v-\varpi(y),z+w-\varpi(z)]_{\Phi+\delta\Phi}.
\end{align*}
This proves the result.
\end{proof}

\begin{prop}
Let  $T:V\rightarrow  \mathfrak{g}$  be  a $\Phi$-twisted Rota-Baxter-operator,
for any  1-cochain $\varpi\in\mathcal{C}^1_{\mathrm{3Lei}}(\mathfrak{g},V)$, if  the linear map $(\mathrm{Id}_V-\varpi\circ T):V\rightarrow V$ is invertible, then the map
$T\circ (\mathrm{Id}_V-\varpi\circ T)^{-1}:V\rightarrow  \mathfrak{g}$ is a $(\Phi+\delta\Phi)$-twisted Rota-Baxter-operator.
\end{prop}

\begin{proof}
  Consider the subalgebra $Gr(T)\subset  \mathfrak{g}\ltimes_{H} V$ of the  $\Phi$-twisted semidirect product.
Thus by Proposition \ref{prop:isomorphism of 3-Leibniz algebras}, we get that
$$\Phi_\varpi(Gr(T))=\{Tu+u-\varpi(Tu)~|~u\in V\}\subset  \mathfrak{g}\ltimes_{\Phi+\delta \varpi}V$$
is a subalgebra. Since the map $(\mathrm{Id}_V-\varpi\circ T):V\rightarrow V$ is invertible, we have $\Phi_\varpi(Gr(T))$
is the graph of the linear map $T\circ(\mathrm{Id}_V-\varpi\circ T)^{-1}$.  Hence by Proposition  \ref{prop:graph}, the map
$T\circ (\mathrm{Id}_V-\varpi\circ T)^{-1}$
is a $(\Phi+\delta\Phi)$-twisted Rota-Baxter-operator.
\end{proof}

Next, we give a construction of a new  $\Phi$-twisted Rota-Baxter-operator out of an old one and a suitable 1-
cocycle.
Let  $T:V\rightarrow  \mathfrak{g}$  be  a $\Phi$-twisted Rota-Baxter-operator.  Suppose $\varpi\in\mathcal{C}^1_{\mathrm{3Lei}}(\mathfrak{g},V)$ is a
1-cocycle. Then $\varpi$ is
said to be $T$-admissible if the linear map $(\mathrm{Id}_V+\varpi\circ T):V\rightarrow V$ is invertible. With this
notation, we have the following.

\begin{prop}
Let $\varpi\in\mathcal{C}^1_{\mathrm{3Lei}}(\mathfrak{g},V)$ be a
$T$-admissible 1-cocycle. Then  the map
$T\circ (\mathrm{Id}_V+\varpi\circ T)^{-1}:V\rightarrow  \mathfrak{g}$ is a  $\Phi$-twisted Rota-Baxter-operator.
\end{prop}

\begin{proof}
Consider the deformed subspace
$$\tau_\varpi(Gr(T))=\{Tu+u+(\varpi\circ T)u~|~u\in V\}\subset  \mathfrak{g}\ltimes_{\Phi}V.$$
Since $\varpi$ is a 1-cocycle, $\tau_\varpi(Gr(T))\subset  \mathfrak{g}\ltimes_{\Phi}V$ turns out to be a subalgebra. Further,
the map $(\mathrm{Id}_V+\varpi\circ T)$ is invertible implies that $\tau_\varpi(Gr(T))$ is the graph of the map
$T\circ (\mathrm{Id}_V+\varpi\circ T)^{-1}$.   Then it follows from
Proposition  \ref{prop:graph} that $T\circ (\mathrm{Id}_V+\varpi\circ T)^{-1}:V\rightarrow  \mathfrak{g}$ is a  $\Phi$-twisted Rota-Baxter-operator.
\end{proof}

The $\Phi$-twisted Rota-Baxter-operator in the above proposition is called the gauge transformation of $T$ associated
with $\varpi$. We denote this $\Phi$-twisted Rota-Baxter-operator simply by $T_\varpi$.

\begin{prop}

Let  $T:V\rightarrow  \mathfrak{g}$  be  a $\Phi$-twisted Rota-Baxter operator and  $\varpi\in\mathcal{C}^1_{\mathrm{3Lei}}(\mathfrak{g},V)$  be a
$T$-admissible 1-cocycle.  Then the  3-Leibniz algebra  structures on $V$ induced from the $\Phi$-twisted Rota-Baxter operators $T$ and $T_\varpi$ are isomorphic.
\end{prop}

\begin{proof}
  Consider the linear isomorphism $(\mathrm{Id}_V+\varpi\circ T): V\rightarrow V$. Moreover, for any
$u,v, w\in V$, we have
\begin{align*}
&[(\mathrm{Id}_V+\varpi\circ T)u, (\mathrm{Id}_V+\varpi\circ T)v,(\mathrm{Id}_V+\varpi\circ T)w]_{T_\varpi}\\
&=\rho^l(Tu,Tv,(\mathrm{Id}_V+\varpi\circ T)w)+\rho^m(Tu,(\mathrm{Id}_V+\varpi\circ T)v,Tw)+\rho^r((\mathrm{Id}_V+\varpi\circ T)u, Tv,Tw)\\
&+\Phi(Tu,Tv,Tw)\\
&=\rho^l(Tu,Tv,w)+\rho^m(Tu,v,Tw)+\rho^r(u, Tv,Tw)+\Phi(Tu,Tv,Tw)+\rho^l(Tu,Tv,(\varpi\circ T)w)\\
&+\rho^m(Tu,(\varpi\circ T)v,Tw)+\rho^r((\varpi\circ T)u, Tv,Tw)\\
&=[u,v,w]_T+\varpi([Tu,Tv,Tw]_\mathfrak{g})\\
&=[u,v,w]_T+\varpi\circ T([u,v,w]_T)\\
&=(\mathrm{Id}_V+\varpi\circ T)[u,v,w]_T.
\end{align*}
Thus $(\mathrm{Id}_V+\varpi\circ T)$ is an isomorphism of 3-Leibniz algebras from $(V,[\cdot,\cdot,\cdot]_T)$ to
$(V,[\cdot,\cdot,\cdot]_{T_\varpi})$.
\end{proof}

\section{Maurer-Cartan characterizations of twisted Rota-Baxter operators }\label{sec: Maurer-Cartan}
\def\theequation{\arabic{section}.\arabic{equation}}
\setcounter{equation} {0}

 In this section, we construct an $L_\infty$-algebra whose Maurer-Cartan elements are twisted Rota-Baxter-operators on 3-Leibniz algebras.
Then we obtain the $L_\infty$-algebra that controls deformations of twisted Rota-Baxter-operators on 3-Leibniz algebras.
 This fact can be viewed as certain justification of twisted Rota-Baxter-operators on 3-Leibniz algebras being interesting structures.

Let $\mathfrak{g}$ be a vector space. We consider the graded vector space $\mathfrak{C}_{\mathrm{3leib}}^\ast(\mathfrak{g},\mathfrak{g})=\oplus_{n\geq -1}\mathfrak{C}_{\mathrm{3Leib}}^{n+1}(\mathfrak{g},\mathfrak{g})$, where
$\mathfrak{C}_{\mathrm{3Leib}}^0(\mathfrak{g},\mathfrak{g})=\otimes^2 \mathfrak{g}$ and
 $\mathfrak{C}_{\mathrm{3Leib}}^{n+1}(\mathfrak{g},\mathfrak{g})$
is the set of linear maps $P\in \mathrm{Hom}\big(\underbrace{(\otimes^2 \mathfrak{g})\otimes\cdots\otimes(\otimes^2 \mathfrak{g})}_n\otimes \mathfrak{g},$ $\mathfrak{g}\big)$.

The degree of elements in $\mathfrak{C}_{\mathrm{3Leib}}^{n+1}(\mathfrak{g},\mathfrak{g})$ are defined to be $n$. Define
\begin{align*}
[P,Q]_{\mathrm{3Leib}}=P\circ Q-(-1)^{pq}Q\circ P, \forall P\in \mathfrak{C}_{\mathrm{3Leib}}^{p+1}(\mathfrak{g},\mathfrak{g}),Q\in \mathfrak{C}_{\mathrm{3Leib}}^{q+1}(\mathfrak{g},\mathfrak{g}),
\end{align*}
where $P\circ Q\in \mathfrak{C}_{\mathrm{3Leib}}^{p+q+1}(\mathfrak{g},\mathfrak{g})$ is defined by
\begin{align*}
&(P\circ Q)(\mathfrak{X}_1,\ldots,\mathfrak{X}_{p+q},x)\\
=&\sum_{k=1}^p(-1)^{(k-1)q}\sum_{\sigma\in S(k-1,q)}(-1)^\sigma P\Big(\mathfrak{X}_{\sigma(1)},\ldots,\mathfrak{X}_{\sigma{(k-1)}},Q\big(\mathfrak{X}_{\sigma(k)},\ldots,\mathfrak{X}_{\sigma{(k+q-1)}},x_{(k+q)}\big)\otimes  y_{k+q},\\
&\mathfrak{X}_{k+q+1},\ldots,\mathfrak{X}_{p+q},x\Big)+\sum_{k=1}^p(-1)^{(k-1)q}\sum_{\sigma\in S(k-1,q)}(-1)^\sigma P\Big(\mathfrak{X}_{\sigma(1)},\ldots,\mathfrak{X}_{\sigma{(k-1)}}, \\
&x_{k+q}\otimes Q\big(\mathfrak{X}_{\sigma(k)},\ldots,\mathfrak{X}_{\sigma{(k+q-1)}},y_{k+q}\big),\mathfrak{X}_{k+q+1},\ldots,\mathfrak{X}_{p+q},x\Big)\\
&+\sum_{\sigma\in S(p,q)}(-1)^{pq}(-1)^\sigma P\Big(\mathcal{X}_{\sigma(1)},\ldots,\mathfrak{X}_{\sigma{(p)}},Q\big(\mathfrak{X}_{\sigma(p+1)},\ldots,\mathfrak{X}_{\sigma{(p+q)}},x\big)\Big),
\end{align*}
for $\mathfrak{X}_i=x_i\otimes y_i\in \otimes^2 \mathfrak{g}, i=1,2,\ldots,p+q$ and $x\in \mathfrak{g}$.

\begin{prop}\cite{Voronov}
With the above notations, the pair $(\mathfrak{C}_{\mathrm{3Leib}}^\ast(\mathfrak{g},\mathfrak{g}),[\cdot,\cdot]_{\mathrm{3Leib}})$   is a graded Lie algebra.
\end{prop}

Let $(V;\rho^l,\rho^m,\rho^r)$ be a representation of a 3-Leibniz algebra $(\mathfrak{g}, [\cdot,\cdot,\cdot]_\mathfrak{g})$ and let $\Phi$ be a 2-cocycle in the cohomology
of $\mathfrak{g}$ with coefficients in $V$. For convenience, we use $\mu: \otimes^3\mathfrak{g}\rightarrow\mathfrak{g}$ to indicate the 3-Leibniz bracket on $\mathfrak{g}$.
In the sequel, we will
view $\rho^l$, $\rho^m$, $\rho^r$ and  $\Phi$ as elements in $\mathrm{Hom}(\otimes^2\mathfrak{g}\otimes V,V)$, $\mathrm{Hom}(\mathfrak{g}\otimes V\otimes\mathfrak{g},V)$,
$\mathrm{Hom}(V\otimes\mathfrak{g}\otimes\mathfrak{g},V)$ and $\mathrm{Hom}(\otimes^3\mathfrak{g},V)$  respectively. Then the $\Phi$-twisted semidirect product 3-Leibniz algebra corresponds to
\begin{align*}
& (\mu+\rho^l+\rho^m+\rho^r+\Phi) \big(x+u,y+v,z+w\big)\\
&=[x,y,z]_\mathfrak{g}+\rho^l(x,y,w)+\rho^m(x,v,z)+\rho^r(u,y,z)+\Phi(x,y,z).
\end{align*}

Therefore, we have $[\mu+\rho^l+\rho^m+\rho^r+\Phi,\mu+\rho^l+\rho^m+\rho^r+\Phi]_{\mathrm{3Leib}}=0.$

\begin{lemma} \label{lemma:ETLTS}
Let $T:V\rightarrow  \mathfrak{g}$ be  a $\Phi$-twisted Rota-Baxter operator on a 3-Leibniz algebra $(\mathfrak{g}, [\cdot,\cdot,\cdot]_\mathfrak{g})$ with respect to the representation $(V;\rho^l,\rho^m,\rho^r)$.  Then, we have
\begin{align*}
&[[\mu+\rho^l+\rho^m+\rho^r+\Phi,T]_{\mathrm{3Leib}},T]_{\mathrm{3Leib}}\big(x+u,y+v,z+w\big)\\
=&2\Big([Tu,Tv,z]_\mathfrak{g}+[Tu,y,Tw]_\mathfrak{g}+[x,Tv,Tw]_\mathfrak{g}-T\big(\rho^l(Tu,y,w)+\rho^m(Tu,v,z)+\\
&\rho^l(x,Tv,w)+\rho^r(u,Tv,z)+\rho^m(x,v,Tw)+\rho^r(u,y,Tw)+\Phi(x,Tv,z)+\Phi(Tu,y,z)+\Phi(x,y,Tw)\big)+\\
&\rho^l(Tu,Tv,w)+\rho^m(Tu,v,Tw)+\rho^r(u,Tv,Tw)+\Phi(Tu,Tv,z)+\Phi(x,Tv,Tw)+\Phi(Tu,y,Tw)
\Big),\\
&[[[\mu+\rho^l+\rho^m+\rho^r+\Phi,T]_{\mathrm{3Leib}},T]_{\mathrm{3Leib}},T]_{\mathrm{3Leib}}\big(x+u,y+v,z+w\big)\\
=&6\Big([Tu,Tv,Tw]_\mathfrak{g}-T\big(\rho^l(Tu,Tv,w)+\rho^m(Tu,v,Tw)+\rho^r(u,Tv,Tw)+\Phi(Tu,Tv,z)+ \\
&\Phi(x,Tv,Tw)+\Phi(Tu,y,Tw)\big)+\Phi(Tu,Tv,Tw)\Big),
\end{align*}
for all $x+u,y+v,z+w\in  \mathfrak{g}\oplus V$.
\end{lemma}

\begin{proof}
For any  $x+u,y+v,z+w\in  \mathfrak{g}\oplus V$ and $T\in \mathfrak{C}^1(V, \mathfrak{g})$, we have
\begin{align*}
&[[\mu+\rho^l+\rho^m+\rho^r+\Phi,T]_{\mathrm{3Leib}},T]_{\mathrm{3Leib}}\big(x+u,y+v,z+w\big)\\
=&[\mu+\rho^l+\rho^m+\rho^r+\Phi,T]_{\mathrm{3Leib}}\big(T(x+u), y+v,z+w\big)+\\
&[\mu+\rho^l+\rho^m+\rho^r+\Phi,T]_{\mathrm{3Leib}}\big(x+u,T(y+v),z+w\big)+\\
&[\mu+\rho^l+\rho^m+\rho^r+\Phi,T]_{\mathrm{3Leib}}\big(x+u,y+v,T(z+w)\big)-\\
&T[\mu+\rho^l+\rho^m+\rho^r+\Phi,T]_{\mathrm{3Leib}}\big(x+u,y+v,z+w\big)\\
=&2\Big([Tu,Tv,z]_\mathfrak{g}+[Tu,y,Tw]_\mathfrak{g}+[x,Tv,Tw]_\mathfrak{g}-T\big(\rho^l(Tu,y,w)+\rho^m(Tu,v,z)+\\
&\rho^l(x,Tv,w)+\rho^r(u,Tv,z)+\rho^m(x,v,Tw)+\rho^r(u,y,Tw)+\Phi(x,Tv,z)+\Phi(Tu,y,z)+\Phi(x,y,Tw)\big)+\\
&\rho^l(Tu,Tv,w)+\rho^m(Tu,v,Tw)+\rho^r(u,Tv,Tw)+\Phi(Tu,Tv,z)+\Phi(x,Tv,Tw)+\Phi(Tu,y,Tw)\Big)
\end{align*}
and
\begin{align*}
&[[[\mu+\rho^l+\rho^m+\rho^r+\Phi,T]_{\mathrm{3Leib}},T]_{\mathrm{3Leib}},T]_{\mathrm{3Leib}}\big(x+u,y+v,z+w\big)\\
=&[[\mu+\rho^l+\rho^m+\rho^r+\Phi,T]_{\mathrm{3Leib}},T]_{\mathrm{3Leib}}\big(T(x+u), y+v,z+w\big)+\\
&[[\mu+\rho^l+\rho^m+\rho^r+\Phi,T]_{\mathrm{3Leib}},T]_{\mathrm{3Leib}}\big(x+u,T(y+v),z+w\big)+\\
&[[\mu+\rho^l+\rho^m+\rho^r+\Phi,T]_{\mathrm{3Leib}},T]_{\mathrm{3Leib}}\big(x+u,y+v,T(z+w)\big)-\\
&T[[\mu+\rho^l+\rho^m+\rho^r+\Phi,T]_{\mathrm{3Leib}},T]_{\mathrm{3Leib}}\big(x+u,y+v,z+w\big)\\
=&6\Big([Tu,Tv,Tw]_\mathfrak{g}-T\big(\rho^l(Tu,Tv,w)+\rho^m(Tu,v,Tw)+\rho^r(u,Tv,Tw)+ \\
&\Phi(Tu,Tv,z)+\Phi(x,Tv,Tw)+\Phi(Tu,y,Tw)\big)+\Phi(Tu,Tv,Tw)\Big).
\end{align*}The proof is finished.
\end{proof}

Consider the graded vector space
$$\mathfrak{C}_{\mathrm{3Leib}}^\ast(V,\mathfrak{g})=\oplus_{n\geq 0}\mathfrak{C}_{\mathrm{3Leib}}^{n+1}(V,\mathfrak{g})=\mathrm{Hom}\big(\underbrace{(\otimes^2 V)\otimes\cdots\otimes(\otimes^2 V)}_n\otimes V, \mathfrak{g}\big).$$

Define
\begin{align*}
&l_3: \mathfrak{C}_{\mathrm{3Leib}}^{m+1}(V,\mathfrak{g})\times\mathfrak{C}_{\mathrm{3Leib}}^{n+1}(V,\mathfrak{g})\times\mathfrak{C}_{\mathrm{3Leib}}^{p+1}(V,\mathfrak{g})\rightarrow\mathfrak{C}_{\mathrm{3Leib}}^{m+n+p+1}(V,\mathfrak{g}),\\
&l_4: \mathfrak{C}_{\mathrm{3Leib}}^{m+1}(V,\mathfrak{g})\times\mathfrak{C}_{\mathrm{3Leib}}^{n+1}(V,\mathfrak{g})\times\mathfrak{C}_{\mathrm{3Leib}}^{p+1}(V,\mathfrak{g})\times\mathfrak{C}_{\mathrm{3Leib}}^{q+1}(V,\mathfrak{g})\rightarrow\mathfrak{C}_{\mathrm{3Leib}}^{m+n+p+q+1}(V,\mathfrak{g})
\end{align*}
by
\begin{align*}
&l_3(P,Q,R)=[[[\mu+\rho^l+\rho^m+\rho^r+\Phi,P]_{\mathrm{3Leib}},Q]_{\mathrm{3Leib}},R]_{\mathrm{3Leib}},\\
&l_4(P,Q,R,S)=[[[[\mu+\rho^l+\rho^m+\rho^r+\Phi,P]_{\mathrm{3Leib}},Q]_{\mathrm{3Leib}},R]_{\mathrm{3Leib}},S]_{\mathrm{3Leib}}.
\end{align*}

In the sequel, we use Voronov's derived brackets theory \cite{Voronov} to construct explicit $L_\infty$-algebras

 \begin{theorem} \label{theorem:Maurer-Cartan equation}
Let $(V;\rho^l,\rho^m,\rho^r)$ be a representation of a 3-Leibniz algebra $(\mathfrak{g}, [\cdot,\cdot,\cdot]_\mathfrak{g})$ and let $\Phi$ be a 2-cocycle in the cohomology
of $\mathfrak{g}$ with coefficients in $V$. Then the graded vector space $\mathfrak{C}_{\mathrm{3Leib}}^\ast(V,\mathfrak{g})$ is an $L_\infty$-algebra with
\begin{align*}
&l_1=l_2=0,~~l_3(\cdot,\cdot,\cdot),~~l_4(\cdot,\cdot,\cdot,\cdot), ~~l_n=0~ \text{for}~ n\geq 5.
\end{align*}
A linear map $T:V\rightarrow  \mathfrak{g}$ is  a $\Phi$-twisted Rota-Baxter operator if and only if $T$ is a solution of the
Maurer-Cartan equation of the $L_\infty$-algebra $(\mathfrak{C}_{\mathrm{3leib}}^\ast(V,\mathfrak{g}),l_3,l_4)$, i.e
\begin{align}
& \frac{1}{3!}l_3(T,T,T)+\frac{1}{4!}l_4(T,T,T,T)=0. \label{4.1}
\end{align}
\end{theorem}
\begin{proof}
By \cite{Voronov},  $(\mathfrak{C}_{\mathrm{3Leib}}^\ast(V,\mathfrak{g}),l_3,l_4)$ is an $L_\infty$-algebra.
 For the second part, by Lemma \ref{lemma:ETLTS}, we have that for any
$T\in\mathfrak{C}_{\mathrm{3Leib}}^1(V,\mathfrak{g})$,
\begin{align}
&l_3(T,T,T)(u,v,w)=6\big([Tu,Tv,Tw]_\mathfrak{g}-T(\rho^l(Tu,Tv,w)+\rho^m(Tu,v,Tw)+\rho^r(u,Tv,Tw))\big) \label{4.2}
\end{align}
and
\begin{align}
&l_4(T,T,T,T)(u,v,w)=-24T\Phi(Tu,Tv,Tw). \label{4.3}
\end{align}
Hence, by Eqs. \eqref{4.2} and \eqref{4.3}, we get
\begin{align*}
& \frac{1}{3!}l_3(T,T,T)(u,v,w)+\frac{1}{4!}l_4(T,T,T,T)(u,v,w)\\
&=[Tu,Tv,Tw]_\mathfrak{g}-T(\rho^l(Tu,Tv,w)+\rho^m(Tu,v,Tw)+\rho^r(u,Tv,Tw)+\Phi(Tu,Tv,Tw)).
\end{align*}
Thus, a linear map $T:V\rightarrow  \mathfrak{g}$ is  a $\Phi$-twisted Rota-Baxter operator if and only if $T$ is a
Maurer-Cartan element of the $L_\infty$-algebra $(\mathfrak{C}_{\mathrm{3Leib}}^\ast(V,\mathfrak{g}),l_3,l_4)$.
\end{proof}

Let $T:V\rightarrow  \mathfrak{g}$ be  a $\Phi$-twisted Rota-Baxter operator on a 3-Leibniz algebra $(\mathfrak{g}, [\cdot,\cdot,\cdot]_\mathfrak{g})$ with respect to the representation $(V;\rho^l,\rho^m,\rho^r)$.
 Since $T$ is a Maurer-Cartan element of the $L_\infty$-algebra $(\mathfrak{C}_{\mathrm{3Leib}}^\ast(V,\mathfrak{g}),l_3,l_4)$,
 then by \cite{Getzler}, $\mathfrak{C}_{\mathrm{3Leib}}^\ast(V, \mathfrak{g})$ carries a
twisted $L_\infty$-algebra structure as following:
\begin{align}
l_1^T(P)=&\frac{1}{2}l_3(T,T,P)+\frac{1}{6}l_4(T,T,T,P), \label{4.4}\\
l_2^T(P,Q)=&l_3(T,P,Q)+\frac{1}{2}l_4(T,T,P,Q),\nonumber\\
l_3^T(P,Q,R)=&l_3(P,Q,R)+l_4(T,P,Q,R),\nonumber\\
l_4^T(P,Q,R,S)=&l_4(P,Q,R,S),\nonumber\\
l_k^T=&0, k\geq 5,\nonumber
\end{align}
where $P\in  \mathfrak{C}_{\mathrm{3Leib}}^{p+1}(V, \mathfrak{g})$, $Q\in  \mathfrak{C}_{\mathrm{3Leib}}^{q+1}(V, \mathfrak{g})$,  $R\in  \mathfrak{C}_{\mathrm{3Leib}}^{r+1}(V, \mathfrak{g})$ and
$S\in  \mathfrak{C}_{\mathrm{3Leib}}^{s+1}(V, \mathfrak{g})$.

Actually, the above twisted $L_\infty$-algebra controls deformations of $\Phi$-twisted Rota-Baxter operators on 3-Leibniz algebras.

\begin{theorem} \label{theorem:MCD}
Let $T:V\rightarrow  \mathfrak{g}$ be  a $\Phi$-twisted Rota-Baxter operator on a 3-Leibniz algebra $(\mathfrak{g}, [\cdot,\cdot,\cdot]_\mathfrak{g})$ with respect to the representation $(V;\rho^l,\rho^m,\rho^r)$.
Then   a linear map $\widetilde{T}:V\rightarrow  \mathfrak{g},$
$T+\widetilde{T}$ is also a  $\Phi$-twisted Rota-Baxter operator
if and only if $\widetilde{T}$ is a Maurer-Cartan element of the twisted $L_\infty$-algebra $(\mathfrak{C}_{\mathrm{3Leib}}^\ast(V, \mathfrak{g}),l_1^T,l_2^T,l_3^T,l_4^T)$,
that is,
$\widetilde{T}$
satisfies the Maurer-Cartan equation:
$$l_1^T(\widetilde{T})+\frac{1}{2!}l_2^T(\widetilde{T},\widetilde{T})+\frac{1}{3!}l_3^T(\widetilde{T},\widetilde{T},\widetilde{T})+\frac{1}{4!}l_4^T(\widetilde{T},\widetilde{T},\widetilde{T},\widetilde{T})=0.$$
\end{theorem}

\begin{proof}
By Theorem \ref{theorem:Maurer-Cartan equation}, $T+\widetilde{T}$ is a $\Phi$-twisted Rota-Baxter operator if and only if
\begin{align*}
& \frac{1}{3!}l_3(T+\widetilde{T},T+\widetilde{T},T+\widetilde{T})+\frac{1}{4!}l_4(T+\widetilde{T},T+\widetilde{T},T+\widetilde{T},T+\widetilde{T})=0.
\end{align*}
Applying Eq.  \eqref{4.1}, the above equation is equivalent to
$$l_1^T(\widetilde{T})+\frac{1}{2!}l_2^T(\widetilde{T},\widetilde{T})+\frac{1}{3!}l_3^T(\widetilde{T},\widetilde{T},\widetilde{T})+\frac{1}{4!}l_4^T(\widetilde{T},\widetilde{T},\widetilde{T},\widetilde{T})=0,$$
 which implies that $\widetilde{T}$ is a Maurer-Cartan element of the twisted $L_\infty$-algebra $(\mathfrak{C}_{\mathrm{3Leib}}^\ast(V, \mathfrak{g}),$ $l_1^T,l_2^T,l_3^T,l_4^T)$.
\end{proof}

\section{ Cohomology of  twisted Rota-Baxter operators }\label{sec: Cohomology}
\def\theequation{\arabic{section}.\arabic{equation}}
\setcounter{equation} {0}

In this section, we construct a representation of the 3-Leibniz algebra $(V,[\cdot,\cdot,\cdot]_T)$
on the vector space $\mathfrak{g}$, and define the cohomology of twisted Rota-Baxter operators
on 3-Leibniz algebras.

\begin{prop}
Let $T:V\rightarrow  \mathfrak{g}$ be  a $\Phi$-twisted Rota-Baxter operator on a 3-Leibniz algebra $(\mathfrak{g}, [\cdot,\cdot,\cdot]_\mathfrak{g})$ with respect to the representation $(V;\rho^l,\rho^m,\rho^r)$.
Define $\rho_T^l: V\otimes V\otimes \mathfrak{g}\rightarrow \mathfrak{g},\rho_T^m:V\otimes \mathfrak{g}\otimes V\rightarrow \mathfrak{g}, \rho_T^r: \mathfrak{g}\otimes V\otimes V\rightarrow \mathfrak{g}$ by
\begin{align}
&\left\{ \begin{array}{lll}
\rho_T^l(u,v,x):=[Tu,Tv,x]_\mathfrak{g}-T(\rho^m(Tu,v,x)+\rho^r(u,Tv,x)+\Phi(Tu,Tv,x)),\\
\rho_T^m(u,x,v):=[Tu,x,Tv]_\mathfrak{g}-T(\rho^l(Tu,x,v)+\rho^r(u,x,Tv)+\Phi(Tu,x,Tv)),\\
\rho_T^r(x,u,v):=[x,Tu,Tv]_\mathfrak{g}-T(\rho^l(x,Tu,v)+\rho^m(x,u,Tv)+\Phi(x,Tu,Tv)),
 \end{array}  \right.\label{5.1}
\end{align}
for all $u,v\in \mathfrak{g}$ and $u\in V.$
then $(\mathfrak{g};\rho_T^l,\rho_T^m,\rho_T^r)$ is a representation of the 3-Leibniz algebra  $(V,[\cdot,\cdot,\cdot]_T)$ on the vector space $\mathfrak{g}$.
\end{prop}

\begin{proof}
For all $x\in \mathfrak{g}$, $u_1,u_2,u_3,u_4\in V,$  by Eqs. \eqref{2.1}, \eqref{3.1}  and \eqref{3.2}, we have
\begin{align*}
&\rho_T^l(u_1,u_2,\rho_T^l(u_3,u_4,x))-\rho_T^l([u_1,u_2,u_3]_{T}, u_4, x)- \rho_T^l(u_3,[u_1,u_2,u_4]_{T},x)-\rho_T^l(u_3,u_4,\rho_T^l(u_1,u_2,x)) \\
=&[Tu_1,Tu_2,[Tu_3,Tu_4,x]_\mathfrak{g}-T\rho^m(Tu_3,u_4,x)-T\rho^r(u_3,Tu_4,x)-T\Phi(Tu_3,Tu_4,x)]_\mathfrak{g}-\\
&T\rho^m(Tu_1,u_2,[Tu_3,Tu_4,x]_\mathfrak{g}-T\rho^m(Tu_3,u_4,x)-T\rho^r(u_3,Tu_4,x)-T\Phi(Tu_3,Tu_4,x))-\\
&T\rho^r(u_1,Tu_2,[Tu_3,Tu_4,x]_\mathfrak{g}-T\rho^m(Tu_3,u_4,x)-T\rho^r(u_3,Tu_4,x)-T\Phi(Tu_3,Tu_4,x))-\\
&T\Phi(Tu_1,Tu_2,[Tu_3,Tu_4,x]_\mathfrak{g}-T\rho^m(Tu_3,u_4,x)-T\rho^r(u_3,Tu_4,x)-T\Phi(Tu_3,Tu_4,x))-\\
&[T[u_1,u_2,u_3]_{T},Tu_4,x]_\mathfrak{g}+T\rho^m(T[u_1,u_2,u_3]_{T},u_4,x)+T\rho^r([u_1,u_2,u_3]_{T},Tu_4,x)+\\
&T\Phi(T[u_1,u_2,u_3]_{T},Tu_4,x)-[Tu_3,T[u_1,u_2,u_4]_{T},x]_\mathfrak{g}+T\rho^m(Tu_3,[u_1,u_2,u_4]_{T},x)+\\
&T\rho^r(u_3,T[u_1,u_2,u_4]_{T},x)+T\Phi(Tu_3,T[u_1,u_2,u_4]_{T},x)-[Tu_3,Tu_4,[Tu_1,Tu_2,x]_\mathfrak{g}-\\
&T\rho^m(Tu_1,u_2,x)-T\rho^r(u_1,Tu_2,x)-T\Phi(Tu_1,Tu_2,x)]_\mathfrak{g}+T\rho^m(Tu_3,u_4,[Tu_1,Tu_2,x]_\mathfrak{g}-\\
&T\rho^m(Tu_1,u_2,x)-T\rho^r(u_1,Tu_2,x)-T\Phi(Tu_1,Tu_2,x))+\\
&T\rho^r(u_3,Tu_4,[Tu_1,Tu_2,x]_\mathfrak{g}-T\rho^m(Tu_1,u_2,x)-T\rho^r(u_1,Tu_2,x)-T\Phi(Tu_1,Tu_2,x))+\\
&T\Phi(Tu_3,Tu_4,[Tu_1,Tu_2,x]_\mathfrak{g}-T\rho^m(Tu_1,u_2,x)-T\rho^r(u_1,Tu_2,x)-T\Phi(Tu_1,Tu_2,x))\\
=&-[Tu_1,Tu_2,T\rho^m(Tu_3,u_4,x)]_\mathfrak{g}-[Tu_1,Tu_2,T\rho^r(u_3,Tu_4,x)]_\mathfrak{g}-[Tu_1,Tu_2,T\Phi(Tu_3,Tu_4,x)]_\mathfrak{g}-\\
&T\rho^m(Tu_1,u_2,[Tu_3,Tu_4,x]_\mathfrak{g})+T\rho^m(Tu_1,u_2,T\rho^m(Tu_3,u_4,x))+T\rho^m(Tu_1,u_2,T\rho^r(u_3,Tu_4,x))+\\
&T\rho^m(Tu_1,u_2,T\Phi(Tu_3,Tu_4,x))-T\rho^r(u_1,Tu_2,[Tu_3,Tu_4,x]_\mathfrak{g})+T\rho^r(u_1,Tu_2,T\rho^m(Tu_3,u_4,x))+\\
&T\rho^r(u_1,Tu_2,T\rho^r(u_3,Tu_4,x))+T\rho^r(u_1,Tu_2,T\Phi(Tu_3,Tu_4,x))+T\Phi(Tu_1,Tu_2,T\rho^m(Tu_3,u_4,x))+\\
&T\Phi(Tu_1,Tu_2,T\rho^r(u_3,Tu_4,x))+T\Phi(Tu_1,Tu_2,T\Phi(Tu_3,Tu_4,x))+T\rho^m([Tu_1,Tu_2,Tu_3]_{\mathfrak{g}},u_4,x)+\\
&T\rho^r(\rho^l(Tu_1,Tu_2,u_3),Tu_4,x)+T\rho^r(\rho^m(Tu_1,u_2,Tu_3),Tu_4,x)+T\rho^r(\rho^r(u_1,Tu_2,Tu_3),Tu_4,x)+\\
&T\rho^m(Tu_3,\rho^l(Tu_1,Tu_2,u_4),x)+T\rho^m(Tu_3,\rho^m(Tu_1,u_2,Tu_4),x)+T\rho^m(Tu_3,\rho^r(u_1,Tu_2,Tu_4),x)+\\
&T\rho^r(u_3,[Tu_1,Tu_2,Tu_4]_{\mathfrak{g}},x)+[Tu_3,Tu_4,T\rho^m(Tu_1,u_2,x)]_\mathfrak{g}+[Tu_3,Tu_4,T\rho^r(u_1,Tu_2,x)]_\mathfrak{g}+\\
&[Tu_3,Tu_4,T\Phi(Tu_1,Tu_2,x)]_\mathfrak{g}+T\rho^m(Tu_3,u_4,[Tu_1,Tu_2,x]_\mathfrak{g})-T\rho^m(Tu_3,u_4,T\rho^m(Tu_1,u_2,x))-\\
&T\rho^m(Tu_3,u_4,T\rho^r(u_1,Tu_2,x))-T\rho^m(Tu_3,u_4,T\Phi(Tu_1,Tu_2,x))+T\rho^r(u_3,Tu_4,[Tu_1,Tu_2,x]_\mathfrak{g})-\\
&T\rho^r(u_3,Tu_4,T\rho^m(Tu_1,u_2,x))-T\rho^r(u_3,Tu_4,T\rho^r(u_1,Tu_2,x))-T\rho^r(u_3,Tu_4,T\Phi(Tu_1,Tu_2,x))-\\
&T\Phi(Tu_3,Tu_4,T\rho^m(Tu_1,u_2,x))-T\Phi(Tu_3,Tu_4,T\rho^r(u_1,Tu_2,x))-T\Phi(Tu_3,Tu_4,T\Phi(Tu_1,Tu_2,x))+\\
&T\rho^l(Tu_1,Tu_2,\Phi(Tu_3,Tu_4,x))-T\rho^l(Tu_3,Tu_4,\Phi(Tu_1,Tu_2,x))\\
=&-T\rho^m(Tu_1,u_2,[Tu_3,Tu_4,x]_\mathfrak{g})-T\rho^r(u_1,Tu_2,[Tu_3,Tu_4,x]_\mathfrak{g})+T\rho^m([Tu_1,Tu_2,Tu_3]_{\mathfrak{g}},u_4,x)+\\
&T\rho^r(\rho^l(Tu_1,Tu_2,u_3),Tu_4,x)+T\rho^r(\rho^m(Tu_1,u_2,Tu_3),Tu_4,x)+T\rho^r(\rho^r(u_1,Tu_2,Tu_3),Tu_4,x)+\\
&T\rho^m(Tu_3,\rho^l(Tu_1,Tu_2,u_4),x)+T\rho^m(Tu_3,\rho^m(Tu_1,u_2,Tu_4),x)+T\rho^m(Tu_3,\rho^r(u_1,Tu_2,Tu_4),x)+\\
&T\rho^r(u_3,[Tu_1,Tu_2,Tu_4]_{\mathfrak{g}},x)+T\rho^m(Tu_3,u_4,[Tu_1,Tu_2,x]_\mathfrak{g})+T\rho^r(u_3,Tu_4,[Tu_1,Tu_2,x]_\mathfrak{g})+\\
&T\rho^l(T u_3,Tu_4,\rho^m(Tu_1,u_2,x))-T\rho^l(T u_1,Tu_2,\rho^m(Tu_3,u_4,x))-T\rho^l(T u_1,Tu_2,\rho^r(u_3,Tu_4,x))+\\
&T\rho^l(T u_3,Tu_4,\rho^r(u_1,Tu_2,x))\\
=&0,
\end{align*}
which implies that Eq. \eqref{2.2} holds. Similarly, we can deduce that Eqs. \eqref{2.3}-\eqref{2.6} hold. Therefore,
 $(\mathfrak{g};\rho_T^l,\rho_T^m,\rho_T^r)$ is a representation of the 3-Leibniz algebra  $(V,[\cdot,\cdot,\cdot]_T)$.
\end{proof}

Let $\delta_T:\mathcal{C}^n_{\mathrm{3Lei}}(V,\mathfrak{g})\rightarrow\mathcal{C}^{n+1}_{\mathrm{3Lei}}(V,\mathfrak{g}), (n \geq 1)$ be the coboundary operator of the
3-Leibniz algebra $(V,[\cdot,\cdot,\cdot]_T)$  with coefficients in the representation $(\mathfrak{g};\rho_T^l,\rho_T^m,\rho_T^r)$. More precisely,
for all $U_i=u_i\otimes v_i\in \otimes^2 V, 1\leq i\leq n$ and $w\in V$, we have
\begin{align*}
&(\delta_T f)(U_1,U_2, \ldots, U_n,w)\\
=&\sum_{1\leq j<k\leq n}(-1)^jf(U_1,\ldots,\widehat{U_j},\cdots,U_{k-1},u_k \otimes[u_j,v_j,v_k]_T+[u_j,v_j,u_k]_T\otimes v_k,\ldots,U_n,w)+\\
&\sum_{j=1}^n(-1)^jf(U_1,\ldots,\widehat{U_j},\ldots,U_{n},[u_j,v_j,w]_T)+\sum_{j=1}^n(-1)^{j+1}\Big([Tu_j,Tv_j,f(U_1,\ldots,\widehat{U_j},\ldots,U_{n},w)]_\mathfrak{g}-\\
&T\rho^m(Tu_j,v_j,f(U_1,\ldots,\widehat{U_j},\ldots,U_{n},w))-T\rho^r(u_j,Tv_j,f(U_1,\ldots,\widehat{U_j},\ldots,U_{n},w))-\\
&T\Phi(Tu_j,Tv_j,f(U_1,\ldots,\widehat{U_j},\ldots,U_{n},w))\Big)+(-1)^{n+1}\Big([Tu_n,f(U_1,\ldots,U_{n-1},v_n),Tw]_\mathfrak{g}-\\
&T\rho^l(Tu_n,f(U_1,\ldots,U_{n-1},v_n),w)-T\rho^r(u_n,f(U_1,\ldots,U_{n-1},v_n),Tw)-\\
&T\Phi(Tu_n,f(U_1,\ldots,U_{n-1},v_n),Tw)+[f(U_1,\ldots,U_{n-1},u_n),Tv_n,Tw]_\mathfrak{g}-\\
&T\rho^l(f(U_1,\ldots,U_{n-1},u_n),Tv_n,w)-T\rho^m(f(U_1,\ldots,U_{n-1},u_n),v_n,Tw)-\\
&T\Phi(f(U_1,\ldots,U_{n-1},u_n),Tv_n,Tw)\Big).
\end{align*}
It is obvious that $f\in \mathcal{C}^1_{\mathrm{3Lei}}(V,\mathfrak{g})$ is closed if and only if
\begin{align*}
&-f\big(\rho^l(Tu,Tv,w)+\rho^m(Tu,v,Tw)+\rho^r(u,Tv,Tw)+\Phi(Tu,Tv,Tw)\big)+\\
&[Tu,Tv,f(w)]_\mathfrak{g}-T\rho^m(Tu,v,f(w))-T\rho^r(u,Tv,f(w))-T\Phi(Tu,Tv,f(w))+\\
&[Tu,f(v),Tw]_\mathfrak{g}-T\rho^l(Tu,f(v),w)-T\rho^r(u,f(v),Tw)-T\Phi(Tu,f(v),Tw)+\\
&[f(u),Tv,Tw]_\mathfrak{g}-T\rho^l(f(u),Tv,w)-T\rho^m(f(u),v,Tw)-T\Phi(f(u_),Tv,Tw)\big)\\
&=0.
\end{align*}

For any $\mathfrak{A}=a\otimes b\in \mathfrak{g}\otimes\mathfrak{g}$, we define $\wp(\mathfrak{A}):V\rightarrow\mathfrak{g}$ by
\begin{align*}
& \wp(\mathfrak{A})(u)=T\rho^l(a,b,u)+T\Phi(a,b,Tu)-[a,b,Tu]_\mathfrak{g}.
\end{align*}

\begin{prop}
Let $T:V\rightarrow  \mathfrak{g}$ be  a $\Phi$-twisted Rota-Baxter operator on a 3-Leibniz algebra $(\mathfrak{g}, [\cdot,\cdot,\cdot]_\mathfrak{g})$ with respect to the representation $(V;\rho^l,\rho^m,\rho^r)$.
 Then $\wp(\mathfrak{A})$ is a 1-cocycle on the 3-Leibniz algebra $(V,[\cdot,\cdot,\cdot]_T)$  with coefficients in   $(\mathfrak{g};\rho_T^l,\rho_T^m,\rho_T^r)$.
\end{prop}

Now, we give a cohomology of $\Phi$-twisted Rota-Baxter operators on a 3-Leibniz algebra.

\begin{definition}
Let $T:V\rightarrow  \mathfrak{g}$ be  a $\Phi$-twisted Rota-Baxter operator on a 3-Leibniz algebra $(\mathfrak{g}, [\cdot,\cdot,\cdot]_\mathfrak{g})$ with respect to the representation $(V;\rho^l,\rho^m,\rho^r)$.
Define the set of $n$-cochains by
\begin{align*}
\mathfrak{C}_T^n(V,\mathfrak{g})=
&\left\{ \begin{array}{ll}
\mathcal{C}^n_{\mathrm{3Lei}}(V,\mathfrak{g}), ~~n\geq 1,\\
\mathfrak{g}\wedge\mathfrak{g},~~~~~~~ ~~n=0.
 \end{array}  \right.
\end{align*}
define $\partial_T:\mathfrak{C}_T^n(V,\mathfrak{g})\rightarrow\mathfrak{C}_T^{n+1}(V,\mathfrak{g})$ by
\begin{align*}
\partial_T=
&\left\{ \begin{array}{ll}
\delta_T, ~~n\geq 1,\\
\wp,~~~n=0.
 \end{array}  \right.
\end{align*}
Denote the set of $n$-cocycles by $\mathrm{Z}_T^n(V,\mathfrak{g})$ and the set of $n$-coboundaries by $\mathrm{B}_T^n(V,\mathfrak{g})$.
Denote by
$$\mathrm{HH}_T^n(V,\mathfrak{g})=\frac{\mathrm{Z}_T^n(V,\mathfrak{g})}{\mathrm{B}_T^n(V,\mathfrak{g})},~~~ n \geq 1$$
the $n$-th cohomology group which will be taken to be the $n$-th cohomology group for the $\Phi$-twisted Rota-Baxter operator $T$.
\end{definition}

At the end of this  section, we give the relationship between the coboundary
operator $\partial_T$ of the $\Phi$-twisted Rota-Baxter operator  $T$ and the differential $l^T_1$ defined
by  Eq.  \eqref{4.4}  using the Maurer-Cartan element $T$ of the $L_\infty$-algebra $(\mathfrak{C}_{\mathrm{3Leib}}^\ast(V,\mathfrak{g}),l_3,l_4)$.

\begin{theorem}\label{theorem:differential}
Let $T:V\rightarrow  \mathfrak{g}$ be  a $\Phi$-twisted Rota-Baxter operator on a 3-Leibniz algebra $(\mathfrak{g}, [\cdot,\cdot,\cdot]_\mathfrak{g})$ with respect to the representation $(V;\rho^l,\rho^m,\rho^r)$.
Then   wa have
$$ l_1^Tf=(-1)^{n-1}\partial_T f,~~\forall f\in \mathrm{Hom}\big(\underbrace{(\otimes^2 V)\otimes\cdots\otimes(\otimes^2 V)}_{n-1}\otimes V, \mathfrak{g}\big).$$
\end{theorem}

\begin{proof}
For all $x,y,z\in \mathfrak{g}$, $u,v,w\in V$, by Lemma \ref{lemma:ETLTS}, we get
\begin{align*}
&[[\mu+\rho^l+\rho^m+\rho^r+\Phi,T]_{\mathrm{3Leib}},T]_{\mathrm{3Leib}}\big(u,v,w\big)=2\big(\rho^l(Tu,Tv,w)+\rho^m(Tu,v,Tw)+\rho^r(u,Tv,Tw)\big),\\
&[[\mu+\rho^l+\rho^m+\rho^r+\Phi,T]_{\mathrm{3Leib}},T]_{\mathrm{3Leib}}\big(u,v,z\big)=2\big([Tu,Tv,z]_\mathfrak{g}-T(\rho^m(Tu,v,z)+\rho^r(u,Tv,z))\big),\\
&[[\mu+\rho^l+\rho^m+\rho^r+\Phi,T]_{\mathrm{3Leib}},T]_{\mathrm{3Leib}}\big(u,y,w\big)=2\big([Tu,y,Tw]_\mathfrak{g}-T\big(\rho^l(Tu,y,w)+\rho^r(u,y,Tw)\big)\big),\\
&[[\mu+\rho^l+\rho^m+\rho^r+\Phi,T]_{\mathrm{3Leib}},T]_{\mathrm{3Leib}}\big(x,v,w\big)=2\big([x,Tv,Tw]_\mathfrak{g}-T\big(\rho^l(x,Tv,w)+\rho^m(x,v,Tw)\big)\big).
\end{align*}
From this, for all $U_i=u_i\otimes v_i\in \otimes^2 V, 1\leq i\leq n$ and $w\in V$,   we have
\begin{align*}
&l_3(T,T,f)(U_1,U_2, \ldots, U_n,w)\\
=&[[\mu+\rho^l+\rho^m+\rho^r+\Phi,T]_{\mathrm{3Leib}},T]_{\mathrm{3Leib}}(f(U_1,U_2, \ldots, u_n)\otimes v_n,w)+\\
&[[\mu+\rho^l+\rho^m+\rho^r+\Phi,T]_{\mathrm{3Leib}},T]_{\mathrm{3Leib}}(u_n\otimes f(U_1,U_2, \ldots, v_n),w)+\\
&\sum_{i=1}^n(-1)^{n+j}[[\mu+\rho^l+\rho^m+\rho^r+\Phi,T]_{\mathrm{3Leib}},T]_{\mathrm{3Leib}}(U_i,f(U_1, \ldots,\widehat{U_i},\ldots, U_n,w))-\\
&(-1)^{n-1}\sum_{1\leq i<k\leq n}(-1)^{i+1}f(U_1,\cdots,\widehat{U_i},\ldots,U_{k}, [[\mu+\rho^l+\rho^m+\rho^r+\Phi,T]_{\mathrm{3Leib}},T]_{\mathrm{3Leib}}(U_i,u_{k+1})\\
&\otimes v_{k+1}+u_{k+1}\otimes[[\mu+\rho^l+\rho^m+\rho^r+\Phi,T]_{\mathrm{3Leib}},T]_{\mathrm{3Leib}}(U_i,v_{k+1}),\ldots, U_n,w)-\\
&(-1)^{n-1}\sum_{i=1}^{n}(-1)^{i+1}f(U_1,\ldots,\widehat{U_i},\ldots,U_n, [[\mu+\rho^l+\rho^m+\rho^r+\Phi,T]_{\mathrm{3Leib}},T]_{\mathrm{3Leib}}(U_i,w))\\
=&2\Big([f(U_1,U_2, \ldots, u_n), Tv_n,Tw]_\mathfrak{g}-T\big(\rho^l(f(U_1,U_2, \ldots, u_n),Tv_n,w)+\rho^m(f(U_1,U_2, \ldots, u_n),v_n,Tw)\big)+\\
&[Tu_n, f(U_1,U_2, \ldots, v_n),Tw]_\mathfrak{g}-T\big(\rho^l(Tu_n,f(U_1,U_2, \ldots, v_n),w)+\rho^r(u_n,f(U_1,U_2, \ldots, v_n),Tw)\big)+\\
&\sum_{i=1}^n(-1)^{n+j-1}\big([Tu_i,Tv_i,f(U_1, \ldots,\widehat{U_i},\ldots, U_n,w)]_\mathfrak{g}-T(\rho^m(Tu_i,v_i,f(U_1, \ldots,\widehat{U_i},\ldots, U_n,w))+\\
&\rho^r(u,Tv,f(U_1, \ldots,\widehat{U_i},\ldots, U_n,w)))\big)-(-1)^{n-1}\sum_{1\leq i<k\leq n}(-1)^{i+1}f(U_1,\ldots,\widehat{U_i},\ldots,U_{k},\\
&(\rho^l(Tu_i,Tv_i,u_{k+1})+\rho^m(Tu_i,v_i,Tu_{k+1})+\rho^r(u_i,Tv_i,Tu_{k+1}))\otimes v_{k+1}+u_{k+1}\otimes(\rho^l(Tu_i,Tv_i,v_{k+1})+\\
&\rho^m(Tu_i,v_i,Tv_{k+1})+\rho^r(u_i,Tv_i,Tv_{k+1})),\ldots, U_n,w)-\\
 &(-1)^{n-1}\sum_{i=1}^{n}(-1)^{i+1}f(U_1,\ldots,\widehat{U_i},\cdots,U_n, \rho^l(Tu_i,Tv_i,w)+\rho^m(Tu_i,v_i,Tw)+\rho^r(u_i,Tv_i,Tw))\Big).
\end{align*}
Also, by Lemma \ref{lemma:ETLTS}, for any $x,y,z\in \mathfrak{g}$, $u,v,w\in V$,  we get
\begin{align*}
&[[[\mu+\rho^l+\rho^m+\rho^r+\Phi,T]_{\mathrm{3Leib}},T]_{\mathrm{3Leib}},T]_{\mathrm{3Leib}}\big(u,v,w\big)=6\Phi(Tu,Tv,Tw),\\
&[[[\mu+\rho^l+\rho^m+\rho^r+\Phi,T]_{\mathrm{3Leib}},T]_{\mathrm{3Leib}},T]_{\mathrm{3Leib}}\big(u,v,z\big)=-6 T\Phi(Tu,Tv,z),\\
&[[[\mu+\rho^l+\rho^m+\rho^r+\Phi,T]_{\mathrm{3Leib}},T]_{\mathrm{3Leib}},T]_{\mathrm{3Leib}}\big(u,y,w\big)=-6 T\Phi(Tu,y,Tw),\\
&[[[\mu+\rho^l+\rho^m+\rho^r+\Phi,T]_{\mathrm{3Leib}},T]_{\mathrm{3Leib}},T]_{\mathrm{3Leib}}\big(x,v,w\big)=-6T\Phi(x,Tv,Tw).
\end{align*}
So we have
\begin{align*}
&l_4(T,T,T,f)(U_1,U_2, \cdots, U_n,w)\\
=&[[[\mu+\rho^l+\rho^m+\rho^r+\Phi,T]_{\mathrm{3Leib}},T]_{\mathrm{3Leib}},T]_{\mathrm{3Leib}}(f(U_1,U_2, \ldots, u_n)\otimes v_n,w)+\\
&[[[\mu+\rho^l+\rho^m+\rho^r+\Phi,T]_{\mathrm{3Leib}},T]_{\mathrm{3Leib}},T]_{\mathrm{3Leib}}(u_n\otimes f(U_1,U_2, \ldots, v_n),w)+\\
&\sum_{i=1}^n(-1)^{n+j}[[[\mu+\rho^l+\rho^m+\rho^r+\Phi,T]_{\mathrm{3Leib}},T]_{\mathrm{3Leib}},T]_{\mathrm{3Leib}}(U_i,f(U_1, \cdots,\widehat{U_i},\cdots, U_n,w))-\\
&(-1)^{n-1}\sum_{1\leq i<k\leq n}(-1)^{i+1}f(U_1,\ldots,\widehat{U_i},\cdots,U_k, [[[\mu+\rho^l+\rho^m+\rho^r+\Phi,T]_{\mathrm{3Leib}},T]_{\mathrm{3Leib}},\\
&T]_{\mathrm{3Leib}}(U_i,u_{k+1})\otimes v_{k+1}+ u_{k+1}\otimes[[[\mu+\rho^l+\rho^m+\rho^r+\Phi,T]_{\mathrm{3Leib}},T]_{\mathrm{3Leib}},T]_{\mathrm{3Leib}}(U_i,v_{k+1}),\ldots, U_n,w)-\\
&(-1)^{n-1}\sum_{i=1}^n(-1)^{i+1}f(U_1,\cdots,\widehat{U_i},\cdots,U_n, [[[\mu+\rho^l+\rho^m+\rho^r+\Phi,T]_{\mathrm{3Leib}},T]_{\mathrm{3Leib}},T]_{\mathrm{3Leib}}(U_i,w))\\
=&6\Big(-T\Phi(f(U_1,U_2, \ldots, u_n),T v_n,Tw)- T\Phi(Tu_n, f(U_1,U_2, \ldots, v_n),Tw)-\\
& \sum_{i=1}^n(-1)^{n+j}T\Phi(Tu_i,Tv_i,f(U_1, \ldots,\widehat{U_i},\ldots, U_n,w))+(-1)^{n-1}\sum_{1\leq i<k\leq n}(-1)^{i+1}f(U_1,\ldots,\widehat{U_i},\\
&\ldots,U_k,\Phi(Tu_i,T v_i,Tu_{k+1}))\otimes v_{k+1}+u_{k+1}\otimes \Phi(Tu_i,T v_i,Tv_{k+1})),\ldots, U_n,w)+\\
&(-1)^{n-1} \sum_{i=1}^{n}(-1)^{i+1}f\big(U_1,\ldots,\widehat{U_i},\cdots,U_n,\Phi(Tu_i,T v_i,Tw)\big)\Big).
\end{align*}
Thus, we deduce that $ l_1^Tf=\frac{1}{2}l_3(T,T,f)+\frac{1}{6}l_4(T,T,T,f)=(-1)^{n-1}\partial_T f.$
\end{proof}

\section{ Deformations of  twisted Rota-Baxter operators }\label{sec: Deformations}
\def\theequation{\arabic{section}.\arabic{equation}}
\setcounter{equation} {0}

In this section, we study formal deformations and extendability of order $n$ deformations to order $n+1$ deformations using the cohomology theory established in the last
section.

Let $\mathbb{K}[[t]]$ be the ring of power series in one variable $t$. For any  linear space $V$, we
let $V[[t]]$ denote the vector space of formal power series in $t$ with coefficients in $V$. If
in addition,$(\mathfrak{g}, [\cdot,\cdot,\cdot]_\mathfrak{g})$ is a  3-Leibniz algebra  over $\mathbb{K}$, then there is a 3-Leibniz algebra structure
over the ring  $\mathbb{K}[[t]]$  on  $ \mathfrak{g}[[t]]$  given by
\begin{align}
[\sum_{i=0}^{\infty}x_it^i,\sum_{j=0}^{\infty}y_jt^j,\sum_{k=0}^{\infty}z_kt^k]_\mathfrak{g}=\sum_{p=0}^{\infty}\sum_{i+j+k=p}[x_i,y_j,z_k]_\mathfrak{g}t^p,\forall x_i,y_j,z_k\in \mathfrak{g}. \label{6.1}
\end{align}

For any representation $(V;\rho^l,\rho^m,\rho^r)$ of $(\mathfrak{g}, [\cdot,\cdot,\cdot]_\mathfrak{g})$, there is a natural representation of the 3-Leibniz
algebra $ \mathfrak{g}[[t]]$  on the $\mathbb{K}[[t]]$ -module $V[[t]]$, which is given by
\begin{align}
&\left\{ \begin{array}{lll}
\rho^l(\sum_{i=0}^{\infty}x_it^i,\sum_{j=0}^{\infty}y_jt^j,\sum_{k=0}^{\infty}v_kt^k)=\sum_{p=0}^{\infty}\sum_{i+j+k=p}\rho^l(x_i,y_j,v_k)t^p,\\
\rho^m(\sum_{i=0}^{\infty}x_it^i,\sum_{k=0}^{\infty}v_kt^k,\sum_{j=0}^{\infty}y_jt^j)=\sum_{p=0}^{\infty}\sum_{i+j+k=p}\rho^m(x_i,v_k,y_j)t^p,\\
\rho^r(\sum_{k=0}^{\infty}v_kt^k,\sum_{i=0}^{\infty}x_it^i,\sum_{j=0}^{\infty}y_jt^j)=\sum_{p=0}^{\infty}\sum_{i+j+k=p}\rho^r(v_k, x_i,y_j)t^p,
 \end{array}  \right.\label{6.2}
\end{align}
for all $x_i,y_j\in \mathfrak{g},z_k\in V.$

Let $T:V\rightarrow  \mathfrak{g}$ be  a $\Phi$-twisted Rota-Baxter operator on a 3-Leibniz algebra $(\mathfrak{g}, [\cdot,\cdot,\cdot]_\mathfrak{g})$ with respect to the representation $(V;\rho^l,\rho^m,\rho^r)$.
Consider a power series
\begin{align}
&T_t=\sum_{i=0}^{\infty}  \mathcal{T}_it^i,~~~\forall \mathcal{T}_i\in \mathrm{Hom}(V,\mathfrak{g}), \label{6.3}
\end{align}
that is, $T_t\in  \mathrm{Hom}(V,\mathfrak{g})[[t]] = \mathrm{Hom}(V,\mathfrak{g}[[t]])$. Extend it to be a $\mathbb{K}[[t]]$-module map
from $V[[t]]$ to $ \mathfrak{g}[[t]]$ which is still denoted by $T_t$.

\begin{definition}
If $T_t=\sum_{i=0}^{\infty}  \mathcal{T}_it^i$ with $\mathcal{T}_0=T$ satisfies
\begin{align}
& [T_tu,T_tv,T_tw]_\mathfrak{g}=T_t(\rho^l(T_tu,T_tv,w)+\rho^m(T_tu,v,T_tw)+\rho^r(u,T_tv,T_tw)+\Phi(T_tu,T_tv,T_tw)),  \label{6.4}
\end{align}
we say that $T_t$ is a formal deformation of the  $\Phi$-twisted Rota-Baxter operator $T$.
\end{definition}

Recall that a formal deformation of a  3-Leibniz algebra $(\mathfrak{g}, [\cdot,\cdot,\cdot]_\mathfrak{g})$ is a power series $\mu_t=\sum_{i=0}^{\infty}  \mu_it^i$,
 where $\mu_i\in \mathrm{Hom}(\mathfrak{g},\mathfrak{g})$ such that $\mu_0=[\cdot,\cdot,\cdot]_\mathfrak{g}$
and $\mu_t$ defines a  3-Leibniz algebra structure   on  $ \mathfrak{g}[[t]]$.

Based on the relationship between relative $\Phi$-twisted Rota-Baxter operators and  3-Leibniz algebras,
we have that follow conclusions.

\begin{prop}
If $T_t=\sum_{i=0}^{\infty}  \mathcal{T}_it^i$ is a formal deformation of a $\Phi$-twisted Rota-Baxter operator $T$ on a  3-Leibniz algebra $(\mathfrak{g}, [\cdot,\cdot,\cdot]_\mathfrak{g})$ with respect to the representation $(V;\rho^l,\rho^m,\rho^r)$, then
$[\cdot,\cdot,\cdot]_{T_t}$ defined by
\begin{align*}
& [u,v,w]_{T_t}= \rho^l(T_tu,T_tv,w)+\rho^m(T_tu,v,T_tw)+\rho^r(u,T_tv,T_tw)+\Phi(T_tu,T_tv,T_tw),
\end{align*}
is a formal deformation of the associated 3-Leibniz algebra $(V,[\cdot,\cdot,\cdot]_{T})$ given in Eq.  \eqref{3.2}.
\end{prop}

Applying Eqs. \eqref{6.1}-\eqref{6.3} to expand Eq. \eqref{6.4} and collecting coefficients of $t^p(p\geq0) $, we see
that Eq.  \eqref{6.4} is equivalent to the system of equations

\begin{align}
\sum_{i+j+k=p}[\mathcal{T}_iu,\mathcal{T}_jv,\mathcal{T}_kw]_\mathfrak{g}=&\sum_{i+j+k=p}\mathcal{T}_i(\rho^l(\mathcal{T}_ju,\mathcal{T}_kv,w)+\rho^m(\mathcal{T}_ju,v,\mathcal{T}_kw)+\rho^r(u,\mathcal{T}_jv,\mathcal{T}_kw))+\label{6.5}\\
&\sum_{i+j+k+l=p}\mathcal{T}_i\Phi(\mathcal{T}_ju,\mathcal{T}_kv,\mathcal{T}_lw),  ~~~~\forall u,v,w\in V.\nonumber
\end{align}

\begin{prop} \label{prop:1-cocycle1}
Let $T_t=\sum_{i=0}^{\infty}  \mathcal{T}_it^i$ be a formal deformation of a $\Phi$-twisted Rota-Baxter operator $T$ on a 3-Leibniz algebra $(\mathfrak{g}, [\cdot,\cdot,\cdot]_\mathfrak{g})$ with respect to the representation $(V;\rho^l,\rho^m,\rho^r)$. Then
$\mathcal{T}_1$ is a 1-cocycle for the $\Phi$-twisted Rota-Baxter operator $T$, that is, $\partial_T\mathcal{T}_1=0$.
\end{prop}

\begin{proof}
For $p=1$, Eq. \eqref{6.5} is equivalent to
 \begin{align*}
&  [\mathcal{T}_1u,Tv,Tw]_\mathfrak{g}+[Tu,\mathcal{T}_1v,Tw]_\mathfrak{g}+[Tu,Tv,\mathcal{T}_1w]_\mathfrak{g}\\
=& \mathcal{T}_1\big(\rho^l(Tu,Tv,w)+\rho^m(Tu,v,Tw)+\rho^r(u,Tv,Tw)+\Phi(Tu,Tv,Tw)\big)+\\
&T\Big(\rho^l(\mathcal{T}_1u,Tv,w)+\rho^m(\mathcal{T}_1u,v,Tw)+\rho^r(u,\mathcal{T}_1v,Tw)+\rho^l(Tu,\mathcal{T}_1v,w)+\rho^m(Tu,v,\mathcal{T}_1w)+\\
&\rho^r(u,Tv,\mathcal{T}_1w)+\Phi(\mathcal{T}_1u,Tv,Tw)+\Phi(Tu,\mathcal{T}_1v,Tw)+\Phi(Tu,Tv,\mathcal{T}_1w)\Big),
\end{align*}
which implies that $\mathcal{T}_1$ is a 1-cocycle.
\end{proof}

\begin{definition}
Let $T$ be  a $\Phi$-twisted Rota-Baxter operator on a 3-Leibniz algebra $(\mathfrak{g}, [\cdot,\cdot,\cdot]_\mathfrak{g})$ with respect to the representation $(V;\rho^l,\rho^m,\rho^r)$.
 The 1-cocycle $\mathcal{T}_1$   given in Proposition \ref{prop:1-cocycle1}  is
called the infinitesimal of the formal deformation $T_t=\sum_{i=0}^{\infty}  \mathcal{T}_it^i$  of $T$.
\end{definition}

\begin{definition}
Two formal deformations $\widetilde{T}_t=\sum_{i=0}^{\infty}  \widetilde{\mathcal{T}_i}t^i$ and $T_t=\sum_{i=0}^{\infty}  \mathcal{T}_it^i$ of a $\Phi$-twisted Rota-Baxter operator
 $T=\widetilde{\mathcal{T}_0}=\mathcal{T}_0$ on a 3-Leibniz algebra $(\mathfrak{g}, [\cdot,\cdot,\cdot]_\mathfrak{g})$ with respect to a representation $(V;\rho^l,\rho^m,\rho^r)$
  are said to be equivalent if there exist $\mathfrak{A}=a\otimes b\in \mathfrak{g}\otimes\mathfrak{g}$, $\eta_i\in \mathrm{End}(\mathfrak{g})$ and
$\zeta_i\in \mathrm{End}(V)$, $i\geq 2$, such that for
 \begin{align}
&f_t=\mathrm{Id}_\mathfrak{g}+t[a,b,\cdot]_\mathfrak{g}+\sum_{i=2}^{\infty}\eta_i t^i,~~g_t=\mathrm{Id}_V+t\big(\rho^l(a,b,\cdot)+\Phi(a,b,T\cdot)\big)+\sum_{i=2}^{\infty}\zeta_i t^i, \label{6.6}
\end{align}
the following conditions hold:
\begin{align*}
&(i)~~f_t[x,y,z]_\mathfrak{g}=[f_t(x),f_t(y),f_t(z)]_\mathfrak{g},\nonumber\\
&(ii)~~\left\{ \begin{array}{ll}
g_t(\rho^l(x,y,u))={\rho}^{l}(f_t(x),f_t(y),g_t(u)),~g_t(\rho^m(x,u,y))={\rho}^{m}(f_t(x),g_t(u),f_t(y)),\\
g_t(\rho^r(u,x,y))={\rho}^r(g_t(u),f_t(x),f_t(y)),~g_t(\Phi(x,y,z))=\Phi(f_t(x),f_t(y),f_t(z)),
 \end{array}  \right.\nonumber \\
&(iii)~~f_t(\widetilde{T_t}u)=T_tg_t(u),~~~~\forall ~~x,y,z\in \mathfrak{g},u\in V.
\end{align*}
\end{definition}

\begin{theorem}
If two formal deformations of a $\Phi$-twisted Rota-Baxter operator $T$ on a  3-Leibniz algebra $(\mathfrak{g}, [\cdot,\cdot,\cdot]_\mathfrak{g})$ with respect to a representation $(V;\rho^l,\rho^m,\rho^r)$ are equivalent, then their infinitesimals are in the same cohomology class.
\end{theorem}

\begin{proof}
Let $(f_t,g_t)$ be the two maps defined by Eq. \eqref{6.6} which gives an equivalence
between two deformations $\widetilde{T}_t=\sum_{i=0}^{\infty}  \widetilde{\mathcal{T}_i}t^i$ and $T_t=\sum_{i=0}^{\infty}  \mathcal{T}_it^i$ of a $\Phi$-twisted Rota-Baxter operator $T$.
By $f_t(\widetilde{T_t}u)=T_tg_t(u)$, we have
$$ \widetilde{\mathcal{T}_1}u=\mathcal{T}_1u+T\rho^l(a,b,u)+T\Phi(a,b,Tu)-[a,b,Tu]_\mathfrak{g}=\mathcal{T}_1u+ \partial_T(\mathfrak{A})(u)$$
which implies that $\widetilde{\mathcal{T}_1}$ and $\mathcal{T}_1$ are in the same cohomology class.
\end{proof}

Next, we introduce a cohomology class associated to any order $n$ deformation of a $\Phi$-twisted Rota-Baxter operator, and show that an order $n$ deformation of a $\Phi$-twisted Rota-Baxter operator
 is extendible if and only if this cohomology class is trivial. Thus we call this
cohomology class the obstruction class of an order $n$ deformation being extendable.

\begin{definition}
Let $T:V\rightarrow  \mathfrak{g}$ be  a $\Phi$-twisted Rota-Baxter operator on a 3-Leibniz algebra $(\mathfrak{g}, [\cdot,\cdot,\cdot]_\mathfrak{g})$ with respect to the representation $(V;\rho^l,\rho^m,\rho^r)$.
If $T_t=\sum_{i=0}^{n}  \mathcal{T}_it^i$ with  $\mathcal{T}_0=T$, $\mathcal{T}_i\in \mathrm{Hom}(V,\mathfrak{g})$,
$i=1,2,\ldots,n$, defines a $\mathbb{K}[t]/(t^{n+1})$-module map
from $V[t]/(t^{n+1})$ to the 3-Leibniz algebra $ \mathfrak{g}[t]/(t^{n+1})$, for all $u,v,w\in V$, satisfying
\begin{align}
& [T_tu,T_tv,T_tw]_\mathfrak{g}=T_t(\rho^l(T_tu,T_tv,w)+\rho^m(T_tu,v,T_tw)+\rho^r(u,T_tv,T_tw)+\Phi(T_tu,T_tv,T_tw)),  \label{6.7}
\end{align}
we say that $T_t$ is an order $n$ deformation of the $\Phi$-twisted Rota-Baxter operator  $T$.
\end{definition}

\begin{definition}
Let $T_t=\sum_{i=0}^{n}  \mathcal{T}_it^i$ be an order $n$ deformation of a $\Phi$-twisted Rota-Baxter operator $T$ on a 3-Leibniz algebra $(\mathfrak{g}, [\cdot,\cdot,\cdot]_\mathfrak{g})$ with respect to the representation $(V;\rho^l,\rho^m,\rho^r)$.
If there
exists a 1-cochain $\mathcal{T}_{n+1}\in \mathrm{Hom}(V,\mathfrak{g})$ such that $\overline{T}_t=T_t+\mathcal{T}_{n+1}t^{n+1}$ is an order $n+1$
deformation of the $\Phi$-twisted Rota-Baxter operator $T$, then we say that $T_t$ is extendable.
\end{definition}

Let $T_t=\sum_{i=0}^{n}  \mathcal{T}_it^i$ be an order $n$ deformation of a $\Phi$-twisted Rota-Baxter operator $T$ on
a  3-Leibniz algebra $(\mathfrak{g}, [\cdot,\cdot,\cdot]_\mathfrak{g})$ with respect to the representation $(V;\rho^l,\rho^m,\rho^r)$.
Define $\mathrm{Ob}_T\in \mathfrak{C}_T^2(V,\mathfrak{g})$
by
  \begin{align}
&\mathrm{Ob}_T(u,v,w)\label{6.8}\\
=&\sum_{\mbox{\tiny$\begin{array}{c}
  i+j+k=n+1\\
   0\leq i,j,k\leq n\end{array}$}}\Big([\mathcal{T}_iu,\mathcal{T}_jv,\mathcal{T}_kw]_\mathfrak{g} -\mathcal{T}_i\big(\rho^l(\mathcal{T}_ju,\mathcal{T}_kv,w)+\rho^m(\mathcal{T}_ju,v,\mathcal{T}_kw)+\rho^r(u,\mathcal{T}_jv,\mathcal{T}_kw)\big)\Big)\nonumber\\
   &-\sum_{\mbox{\tiny$\begin{array}{c}
  i+j+k+l=n+1\\
   0\leq i,j,k,l\leq n\end{array}$}}\mathcal{T}_i\Phi(\mathcal{T}_ju,\mathcal{T}_kv,\mathcal{T}_lw).\nonumber
\end{align}

Using the $L_\infty$-algebra given in Theorem \ref{theorem:Maurer-Cartan equation}, we show that  $\mathrm{Ob}_T$  is a 2-cocycle.

\begin{prop} \label{prop:1-cocycle}
The 2-cochain $\mathrm{Ob}_T$ is a 2-cocycle, that is, $\partial_T\mathrm{Ob}_T=0$.
\end{prop}

\begin{proof}
For any $u,v,w\in V$,   we have
  \begin{align*}
&l_3(\mathcal{T}_i,\mathcal{T}_j,\mathcal{T}_k)(u,v,w)\\
=&[[[\mu+\rho^l+\rho^m+\rho^r+\Phi,\mathcal{T}_i]_{\mathrm{3Leib}},\mathcal{T}_j]_{\mathrm{3Leib}},\mathcal{T}_k]_{\mathrm{3Leib}}(u,v,w)\\
=&[\mu+\rho^l+\rho^m+\rho^r+\Phi,\mathcal{T}_i]_{\mathrm{3Leib}}(\mathcal{T}_k u,\mathcal{T}_jv,w)+[\mu+\rho^l+\rho^m+\rho^r+\Phi,\mathcal{T}_i]_{\mathrm{3Leib}}(\mathcal{T}_k u,v,\mathcal{T}_jw)-\\
&\mathcal{T}_j[\mu+\rho^l+\rho^m+\rho^r+\Phi,\mathcal{T}_i]_{\mathrm{3Leib}}(\mathcal{T}_k u,v,w)+[\mu+\rho^l+\rho^m+\rho^r+\Phi,\mathcal{T}_i]_{\mathrm{3Leib}}(\mathcal{T}_j u,\mathcal{T}_k v,w)+\\
&[\mu+\rho^l+\rho^m+\rho^r+\Phi,\mathcal{T}_i]_{\mathrm{3Leib}}( u,\mathcal{T}_k v,\mathcal{T}_j w)-\mathcal{T}_j [\mu+\rho^l+\rho^m+\rho^r+\Phi,\mathcal{T}_i]_{\mathrm{3Leib}}(u,\mathcal{T}_k v,w)+\\
&[\mu+\rho^l+\rho^m+\rho^r+\Phi,\mathcal{T}_i]_{\mathrm{3Leib}}(\mathcal{T}_j u,v,\mathcal{T}_k w)+[\mu+\rho^l+\rho^m+\rho^r+\Phi,\mathcal{T}_i]_{\mathrm{3Leib}}(u,\mathcal{T}_j v,\mathcal{T}_k w)-\\
&\mathcal{T}_j [\mu+\rho^l+\rho^m+\rho^r+\Phi,\mathcal{T}_i]_{\mathrm{3Leib}}(u,v,\mathcal{T}_k w)-\mathcal{T}_k[\mu+\rho^l+\rho^m+\rho^r+\Phi,\mathcal{T}_i]_{\mathrm{3Leib}}(\mathcal{T}_j u,v,w)-\\
&\mathcal{T}_k[\mu+\rho^l+\rho^m+\rho^r+\Phi,\mathcal{T}_i]_{\mathrm{3Leib}}(u,\mathcal{T}_j v,w)-\mathcal{T}_k[\mu+\rho^l+\rho^m+\rho^r+\Phi,\mathcal{T}_i]_{\mathrm{3Leib}}(u,v,\mathcal{T}_j w)\\
=&[\mathcal{T}_k u,\mathcal{T}_jv,\mathcal{T}_iw]_\mathfrak{g}-\mathcal{T}_i\rho^l(\mathcal{T}_k u,\mathcal{T}_jv,w)+[\mathcal{T}_k u,\mathcal{T}_iv,\mathcal{T}_jw]_\mathfrak{g}-\mathcal{T}_i\rho^m(\mathcal{T}_k u,v,\mathcal{T}_jw)-\mathcal{T}_j\big(\rho^l(\mathcal{T}_k u,\mathcal{T}_iv,w)+\\
&\rho^m(\mathcal{T}_k u,v,\mathcal{T}_iw)\big)+[\mathcal{T}_j u,\mathcal{T}_kv,\mathcal{T}_iw]_\mathfrak{g}-\mathcal{T}_i\rho^l(\mathcal{T}_j u,\mathcal{T}_kv,w)+[\mathcal{T}_i u,\mathcal{T}_kv,\mathcal{T}_jw]_\mathfrak{g}-\mathcal{T}_i\rho^r(u, \mathcal{T}_k v,\mathcal{T}_jw)-\\
&\mathcal{T}_j\big(\rho^l(\mathcal{T}_i u,\mathcal{T}_kv,w)+\rho^r( u,\mathcal{T}_k v,\mathcal{T}_iw)\big)+[\mathcal{T}_j u,\mathcal{T}_iv,\mathcal{T}_kw]_\mathfrak{g}-\mathcal{T}_i\rho^m(\mathcal{T}_j u,v,\mathcal{T}_kw)+
[\mathcal{T}_i u,\mathcal{T}_jv,\mathcal{T}_kw]_\mathfrak{g}-\\
&\mathcal{T}_i\rho^r(u, \mathcal{T}_j v,\mathcal{T}_kw)-\mathcal{T}_j\big(\rho^m(\mathcal{T}_i u,v,\mathcal{T}_kw)+\rho^r( u,\mathcal{T}_i v,\mathcal{T}_kw)\big)-\mathcal{T}_k\big(\rho^l(\mathcal{T}_j u,\mathcal{T}_iv,w)+\rho^m(\mathcal{T}_j u, v,\mathcal{T}_iw)\big)-\\
&\mathcal{T}_k\big(\rho^l(\mathcal{T}_i u,\mathcal{T}_jv,w)+\rho^r( u, \mathcal{T}_j v,\mathcal{T}_iw)\big)-\mathcal{T}_k\big(\rho^r( u,\mathcal{T}_iv,\mathcal{T}_j w)+\rho^m(\mathcal{T}_i u, v,\mathcal{T}_jw)\big)
\end{align*}
and
\begin{align*}
&l_4(\mathcal{T}_i,\mathcal{T}_j,\mathcal{T}_k,\mathcal{T}_l)(u,v,w)\\
=&[[[[\mu+\rho^l+\rho^m+\rho^r+\Phi,\mathcal{T}_i]_{\mathrm{3Leib}},\mathcal{T}_j]_{\mathrm{3Leib}},\mathcal{T}_k]_{\mathrm{3Leib}},\mathcal{T}_l]_{\mathrm{3Leib}}(u,v,w)\\
=&[[[\mu+\rho^l+\rho^m+\rho^r+\Phi,\mathcal{T}_i]_{\mathrm{3Leib}},\mathcal{T}_j]_{\mathrm{3Leib}},\mathcal{T}_k]_{\mathrm{3Leib}}(\mathcal{T}_l u,v,w)+\\
&[[[\mu+\rho^l+\rho^m+\rho^r+\Phi,\mathcal{T}_i]_{\mathrm{3Leib}},\mathcal{T}_j]_{\mathrm{3Leib}},\mathcal{T}_k]_{\mathrm{3Leib}}(u,\mathcal{T}_l v,w)+\\
&[[[\mu+\rho^l+\rho^m+\rho^r+\Phi,\mathcal{T}_i]_{\mathrm{3Leib}},\mathcal{T}_j]_{\mathrm{3Leib}},\mathcal{T}_k]_{\mathrm{3Leib}}(u,v,\mathcal{T}_l w)-\\
&\mathcal{T}_l[[[\mu+\rho^l+\rho^m+\rho^r+\Phi,\mathcal{T}_i]_{\mathrm{3Leib}},\mathcal{T}_j]_{\mathrm{3Leib}},\mathcal{T}_k]_{\mathrm{3Leib}}(u,v,w)\\
=&-\mathcal{T}_i\Phi(\mathcal{T}_l u,\mathcal{T}_k v,\mathcal{T}_jw)-\mathcal{T}_i\Phi(\mathcal{T}_l u,\mathcal{T}_j v,\mathcal{T}_kw)-\mathcal{T}_i\Phi(\mathcal{T}_k u,\mathcal{T}_l v,\mathcal{T}_jw)
-\mathcal{T}_i\Phi(\mathcal{T}_j u,\mathcal{T}_l v,\mathcal{T}_kw)-\\
&\mathcal{T}_i\Phi(\mathcal{T}_k u,\mathcal{T}_j v,\mathcal{T}_lw)-\mathcal{T}_i\Phi(\mathcal{T}_j u,\mathcal{T}_k v,\mathcal{T}_lw)-\mathcal{T}_j\Phi(\mathcal{T}_l u,\mathcal{T}_k v,\mathcal{T}_iw)
-\mathcal{T}_j\Phi(\mathcal{T}_l u,\mathcal{T}_i v,\mathcal{T}_kw)-\\
&\mathcal{T}_k\Phi(\mathcal{T}_l u,\mathcal{T}_j v,\mathcal{T}_iw)-\mathcal{T}_k\Phi(\mathcal{T}_l u,\mathcal{T}_i v,\mathcal{T}_jw)-\mathcal{T}_j\Phi(\mathcal{T}_k u,\mathcal{T}_l v,\mathcal{T}_iw)
-\mathcal{T}_j\Phi(\mathcal{T}_i u,\mathcal{T}_l v,\mathcal{T}_kw)-\\
&\mathcal{T}_k\Phi(\mathcal{T}_j u,\mathcal{T}_l v,\mathcal{T}_iw)-\mathcal{T}_k\Phi(\mathcal{T}_i u,\mathcal{T}_l v,\mathcal{T}_jw)-\mathcal{T}_j\Phi(\mathcal{T}_k u,\mathcal{T}_i v,\mathcal{T}_lw)
-\mathcal{T}_j\Phi(\mathcal{T}_i u,\mathcal{T}_k v,\mathcal{T}_lw)-\\
&\mathcal{T}_k\Phi(\mathcal{T}_j u,\mathcal{T}_i v,\mathcal{T}_lw)-\mathcal{T}_k\Phi(\mathcal{T}_i u,\mathcal{T}_j v,\mathcal{T}_lw)-\mathcal{T}_l\Phi(\mathcal{T}_k u,\mathcal{T}_j v,\mathcal{T}_iw)
-\mathcal{T}_l\Phi(\mathcal{T}_k u,\mathcal{T}_i v,\mathcal{T}_jw)-\\
&\mathcal{T}_l\Phi(\mathcal{T}_j u,\mathcal{T}_k v,\mathcal{T}_iw)-\mathcal{T}_l\Phi(\mathcal{T}_i u,\mathcal{T}_k v,\mathcal{T}_jw)-\mathcal{T}_l\Phi(\mathcal{T}_j u,\mathcal{T}_i v,\mathcal{T}_kw)
-\mathcal{T}_l\Phi(\mathcal{T}_i u,\mathcal{T}_j v,\mathcal{T}_kw).
\end{align*}
Thus, we deduce that
  \begin{align}
\mathrm{Ob}_T =&\frac{1}{6}\sum_{\mbox{\tiny$\begin{array}{c}
  i+j+k=n+1\\
   0\leq i,j,k\leq n\end{array}$}}l_3(\mathcal{T}_i,\mathcal{T}_j,\mathcal{T}_k)+\frac{1}{24}\sum_{\mbox{\tiny$\begin{array}{c}
  i+j+k+l=n+1\\
   0\leq i,j,k,l\leq n\end{array}$}}l_4(\mathcal{T}_i,\mathcal{T}_j,\mathcal{T}_k,\mathcal{T}_l)\label{6.9}\\
   =&\frac{1}{6}\sum_{\mbox{\tiny$\begin{array}{c}
  i+j+k=n+1\\
   1\leq i,j,k\leq n\end{array}$}}l_3(\mathcal{T}_i,\mathcal{T}_j,\mathcal{T}_k)+\frac{1}{24}\sum_{\mbox{\tiny$\begin{array}{c}
  i+j+k+l=n+1\\
   1\leq i,j,k,l\leq n\end{array}$}}l_4(\mathcal{T}_i,\mathcal{T}_j,\mathcal{T}_k,\mathcal{T}_l)+\nonumber\\
   &\frac{1}{2}\sum_{\mbox{\tiny$\begin{array}{c}
  i+j=n+1\\
   1\leq i,j\leq n\end{array}$}}l_3(T,\mathcal{T}_i,\mathcal{T}_j)+\frac{1}{6}\sum_{\mbox{\tiny$\begin{array}{c}
  i+j+k=n+1\\
   1\leq i,j,k\leq n\end{array}$}}l_4(T,\mathcal{T}_i,\mathcal{T}_j,\mathcal{T}_k).\nonumber
\end{align}
Since $T_t$ is an order $n$ deformation of the $\Phi$-twisted Rota-Baxter operator $T$, for all $0\leq s\leq n$, $u,v,w\in V$ we have
  \begin{align*}
&\sum_{\mbox{\tiny$\begin{array}{c}
  i+j+k=s\\
   0\leq i,j,k\leq s\end{array}$}}\Big([\mathcal{T}_iu,\mathcal{T}_jv,\mathcal{T}_kw]_\mathfrak{g} -\mathcal{T}_i\big(\rho^l(\mathcal{T}_ju,\mathcal{T}_kv,w)+\rho^m(\mathcal{T}_ju,v,\mathcal{T}_kw)+\rho^r(u,\mathcal{T}_jv,\mathcal{T}_kw)\big)\Big)- \\
   &\sum_{\mbox{\tiny$\begin{array}{c}
  i+j+k+l=s\\
   0\leq i,j,k,l\leq s\end{array}$}}\mathcal{T}_i\Phi(\mathcal{T}_ju,\mathcal{T}_kv,\mathcal{T}_lw)=0.
\end{align*}
which is equivalent to
  \begin{align}
-\frac{1}{2}l_3(T,T,\mathcal{T}_s)-\frac{1}{6}l_4(T,T,T,\mathcal{T}_s) =&\frac{1}{6}\sum_{\mbox{\tiny$\begin{array}{c}
  i+j+k=s\\
   0\leq i,j,k\leq s-1\end{array}$}}l_3(\mathcal{T}_i,\mathcal{T}_j,\mathcal{T}_k)+\frac{1}{24}\sum_{\mbox{\tiny$\begin{array}{c}
  i+j+k+l=s\\
   0\leq i,j,k,l\leq s-1\end{array}$}}l_4(\mathcal{T}_i,\mathcal{T}_j,\mathcal{T}_k,\mathcal{T}_l)\label{6.10}\\
   =&\frac{1}{6}\sum_{\mbox{\tiny$\begin{array}{c}
  i+j+k=s\\
   1\leq i,j,k\leq s-1\end{array}$}}l_3(\mathcal{T}_i,\mathcal{T}_j,\mathcal{T}_k)+\frac{1}{24}\sum_{\mbox{\tiny$\begin{array}{c}
  i+j+k+l=s\\
   1\leq i,j,k,l\leq s-1\end{array}$}}l_4(\mathcal{T}_i,\mathcal{T}_j,\mathcal{T}_k,\mathcal{T}_l)+\nonumber\\
   &\frac{1}{2}\sum_{\mbox{\tiny$\begin{array}{c}
  i+j=s\\
   1\leq i,j\leq s-1\end{array}$}}l_3(T,\mathcal{T}_i,\mathcal{T}_j)+\frac{1}{6}\sum_{\mbox{\tiny$\begin{array}{c}
  i+j+k=s\\
   1\leq i,j,k\leq s-1\end{array}$}}l_4(T,\mathcal{T}_i,\mathcal{T}_j,\mathcal{T}_k).\nonumber
\end{align}
By Theorem \ref{theorem:differential}, Eq. \eqref{4.4}, \eqref{6.9} and \eqref{6.10}, we have
  \begin{align*}
  \partial_T\mathrm{Ob}_T=&-\frac{1}{2}l_3(T,T,\mathrm{Ob}_T)-\frac{1}{6}l_4(T,T,T,\mathrm{Ob}_T)\\
  =&-\frac{1}{12}\sum_{\mbox{\tiny$\begin{array}{c}
  i+j+k=n+1\\
   1\leq i,j,k\leq n\end{array}$}}l_3(T,T,l_3(\mathcal{T}_i,\mathcal{T}_j,\mathcal{T}_k))-
   \frac{1}{48}\sum_{\mbox{\tiny$\begin{array}{c}
  i+j+k+l=n+1\\
   1\leq i,j,k,l\leq n\end{array}$}}l_3(T,T,l_4(\mathcal{T}_i,\mathcal{T}_j,\mathcal{T}_k,\mathcal{T}_l))-\\
   &\frac{1}{4}\sum_{\mbox{\tiny$\begin{array}{c}
  i+j=n+1\\
   1\leq i,j\leq n\end{array}$}}l_3(T,T,l_3(T,\mathcal{T}_i,\mathcal{T}_j))-
   \frac{1}{12}\sum_{\mbox{\tiny$\begin{array}{c}
  i+j+k=n+1\\
   1\leq i,j,k\leq n\end{array}$}}l_3(T,T,l_4(T,\mathcal{T}_i,\mathcal{T}_j,\mathcal{T}_k))-\\
   &\frac{1}{36}\sum_{\mbox{\tiny$\begin{array}{c}
  i+j+k=n+1\\
   1\leq i,j,k\leq n\end{array}$}}l_4(T,T,T,l_3(\mathcal{T}_i,\mathcal{T}_j,\mathcal{T}_k))-
   \frac{1}{144}\sum_{\mbox{\tiny$\begin{array}{c}
  i+j+k+l=n+1\\
   1\leq i,j,k,l\leq n\end{array}$}}l_4(T,T,T,l_4(\mathcal{T}_i,\mathcal{T}_j,\mathcal{T}_k,\mathcal{T}_l))-\\
   &\frac{1}{12}\sum_{\mbox{\tiny$\begin{array}{c}
  i+j=n+1\\
   1\leq i,j\leq n\end{array}$}}l_4(T,T,T,l_3(T,\mathcal{T}_i,\mathcal{T}_j))-
   \frac{1}{36}\sum_{\mbox{\tiny$\begin{array}{c}
  i+j+k=n+1\\
   1\leq i,j,k\leq n\end{array}$}}l_4(T,T,T,l_4(T,\mathcal{T}_i,\mathcal{T}_j,\mathcal{T}_k))\\
   =&0.
    \end{align*}
    Thus, we obtain that   $\mathrm{Ob}_T$ is a 2-cocycle.
\end{proof}

\begin{definition}
Let $T_t=\sum_{i=0}^{n}  \mathcal{T}_it^i$ be an order $n$ deformation of a $\Phi$-twisted Rota-Baxter operator $T$ on a 3-Leibniz algebra $(\mathfrak{g}, [\cdot,\cdot,\cdot]_\mathfrak{g})$ with respect to the representation $(V;\rho^l,\rho^m,\rho^r)$.
The
cohomology class $[\mathrm{Ob}_T]\in \mathrm{HH}_T^2(V,\mathfrak{g})$  is called the obstruction class of $T_t$ being extendable.
\end{definition}

\begin{theorem}\label{theorem:differential}
Let $T_t=\sum_{i=0}^{n}  \mathcal{T}_it^i$ be an order $n$ deformation of a $\Phi$-twisted Rota-Baxter operator $T$ on a 3-Leibniz algebra $(\mathfrak{g}, [\cdot,\cdot,\cdot]_\mathfrak{g})$ with respect to the representation $(V;\rho^l,\rho^m,\rho^r)$.
Then
$T_t$ is extendable if and only if the obstruction class $[\mathrm{Ob}_T]$ is trivial.
\end{theorem}

\begin{proof}
Suppose that an order $n$ deformation $T_t$ of the $\Phi$-twisted Rota-Baxter operator $T$
extends to an order $n+1$ deformation. Then Eq. \eqref{6.10} holds for $s=n+1$. Thus, we have
$\mathrm{Ob}_T=-\partial_T \mathcal{T}_{n+1}$, which implies that the obstruction class $[\mathrm{Ob}_T]$ is trivial.

Conversely, if the obstruction class $[\mathrm{Ob}_T]$ is trivial, suppose that $\mathrm{Ob}_T=-\partial_T \mathcal{T}_{n+1}$ for
some 1-cochain $\mathcal{T}_{n+1}\in \mathrm{Hom}(V,\mathfrak{g})$. Set $\overline{T}_t=T_t+\mathcal{T}_{n+1}t^{n+1}$.
Then $\overline{T}_t$ satisfies Eq. \eqref{6.10} for
$0\leq s\leq n+1$. So $\overline{T}_t$ is an order $n+1$ deformation, which means that $T_t$ is extendable.
\end{proof}

\begin{coro}
$T:V\rightarrow  \mathfrak{g}$ be  a $\Phi$-twisted Rota-Baxter operator on a 3-Leibniz algebra $(\mathfrak{g}, [\cdot,\cdot,\cdot]_\mathfrak{g})$ with respect to a representation $(V;\rho^l,\rho^m,\rho^r)$.
If $\mathrm{HH}_T^2(V,\mathfrak{g})=0$, then every 1-cocycle in $\mathrm{Z}_T^1(V,\mathfrak{g})$
is the infinitesimal of some formal deformation of the $\Phi$-twisted Rota-Baxter operator $T$.
\end{coro}

\section{ NS-3-Leibniz algebras }\label{sec: NS-3-Leibniz algebras}
\def\theequation{\arabic{section}.\arabic{equation}}
\setcounter{equation} {0}

In this section, we introduce the notion of an  NS-3-Leibniz algebra  as the underlying structure of twisted Rota-Baxter operators.
We study some properties of NS-3-Leibniz algebras and give some examples. Further study on
NS-3-Leibniz algebras is postponed to a forthcoming paper.

\begin{definition}
An NS-3-Leibniz algebra is a vector space $\mathfrak{g}$ endowed with   four  trilinear maps $[\cdot,\cdot,\cdot]_{\triangleleft},[\cdot,\cdot,\cdot]_{\triangleright},[\cdot,\cdot,\cdot]_{\triangle},[\cdot,\cdot,\cdot]_{\diamond}:\mathfrak{g}\times\mathfrak{g}\times\mathfrak{g}\rightarrow\mathfrak{g}$
satisfying for all $ x, y, z, a, b\in \mathfrak{g}$,
\begin{align}
 [a,b,[x,y,z]_\triangleright]_{\triangleright}=&[[a, b, x]_{\star},y,z]_{\triangleright}+ [x,  [a, b, y]_{\star},z]_{\triangleright}+ [x,y,[a, b, z]_{\triangleright}]_{\triangleright},\label{7.1}\\
 [a,b,[x,y,z]_\star]_{\triangleleft}=&[[a, b, x]_{\triangleleft},y,z]_{\triangleleft}+ [x,  [a, b, y]_{\triangleleft},z]_{\triangle}+ [x,y,[a, b, z]_{\triangleleft}]_{\triangleright},\label{7.2}\\
 [a,b,[x,y,z]_\star]_{\triangle}=&[[a, b, x]_{\triangle},y,z]_{\triangleleft}+ [x,  [a, b, y]_{\triangle},z]_{\triangle}+ [x,y,[a, b, z]_{\triangle}]_{\triangleright},\label{7.3}\\
  [a,b,[x,y,z]_\triangleleft]_{\triangleright}=&[[a, b, x]_{\triangleright},y,z]_{\triangleleft}+ [x,  [a, b, y]_{\star},z]_{\triangleleft}+ [x,y,[a, b, z]_{\star}]_{\triangleleft},\label{7.4}\\
    [a,b,[x,y,z]_\diamond]_{\triangleright}=&[[a, b, x]_{\triangleright},y,z]_{\triangleleft}+ [x,  [a, b, y]_{\star},z]_{\triangleleft}+ [x,y,[a, b, z]_{\star}]_{\triangleleft},\label{7.5}\\
[a,b,[x, y, z]_{\star}]_{\diamond}+[a,b,[x, y, z]_{\diamond}]_{\triangleright}=&[x,y,[a,b,z]_\diamond]_{\triangleright}+[x,y,[a,b,z]_\star]_{\diamond}+[x,  [a,b,y]_{\diamond},z]_{\triangle}+\label{7.6}\\
&[x,[a,b, y]_{\star},z]_{\diamond}+[[a,b,x]_{\star},y, z]_{\diamond}+[[a,b,x]_{\diamond},y, z]_{\triangleleft},\nonumber
\end{align}
where $[x,y,z]_\star=[x,y,z]_\triangleleft+[x,y,z]_\triangleright+[x,y,z]_{\bigtriangleup}+[x,y,z]_\diamond.$
\end{definition}

\begin{remark}
If the trilinear operations $[\cdot,\cdot,\cdot]_{\triangleleft}, [\cdot,\cdot,\cdot]_{\triangle}$ and $[\cdot,\cdot,\cdot]_{\diamond}$ in the above definition are trivial,
 it follows that
$(\mathfrak{g},[\cdot,\cdot,\cdot]_{\triangleright})$ is a 3-Leibniz algebra.
On the other hand, if $[\cdot,\cdot,\cdot]_{\diamond}$ is trivial then $(\mathfrak{g},[\cdot,\cdot,\cdot]_{\triangleleft},[\cdot,\cdot,\cdot]_{\triangleright}, [\cdot,\cdot,\cdot]_{\triangle})$  becomes a 3-pre-Leibniz algebra. Thus,
NS-3-Leibniz algebras are a generalization of both 3-Leibniz algebras and 3-pre-Leibniz algebras (For more details
about 3-pre-Leibniz algebras, see \cite{Hu24}).
\end{remark}

\begin{remark}
In \cite{Hou21,Chtioui}, the authors introduced the notion of an NS-3-Lie algebra which is the   underlying structure of a $\Phi$-twisted Rota-Baxter operator  (Generalized Reynolds operator) on a 3-Lie algebra.
In our definition of an NS-3-Leibniz algebras
$(\mathfrak{g},[\cdot,\cdot,\cdot]_{\triangleright},[\cdot,\cdot,\cdot]_{\triangleleft}, [\cdot,\cdot,\cdot]_{\triangle}, [\cdot,\cdot,\cdot]_{\diamond})$, if the ternary product $[\cdot,\cdot,\cdot]_\triangleright:\mathfrak{g}\times\mathfrak{g}\times\mathfrak{g}\rightarrow\mathfrak{g}$
satisfies $[x,y,z]_\triangleright=-[y,x,z]_\triangleright$, $[x,y,z]_\bigtriangleup=[z,x,y]_\triangleright$ and $[x,y,z]_\triangleleft=[y,z,x]_\triangleright$ for $x,y,z\in\mathfrak{g}$, then $(\mathfrak{g},[\cdot,\cdot,\cdot]_{\triangleright}, [\cdot,\cdot,\cdot]_{\diamond})$
becomes an NS-3-Lie algebra. So the NS-3-Leibniz algebras are natural generalizations of NS-3-Lie algebras.
\end{remark}

In the following, we show that NS-3-Leibniz algebras split 3-Leibniz algebras.

\begin{prop} \label{prop:sub-adjacent 3-Leibniz}
Let $(\mathfrak{g},[\cdot,\cdot,\cdot]_{\triangleright},  [\cdot,\cdot,\cdot]_{\triangleleft}, [\cdot,\cdot,\cdot]_{\triangle}, [\cdot,\cdot,\cdot]_{\diamond})$  be an NS-3-Leibniz algebra. Then $(\mathfrak{g},[\cdot,\cdot,\cdot]_{\star})$  is a 3-Leibniz algebra
which is called the sub-adjacent 3-Leibniz  algebra of $(\mathfrak{g},[\cdot,\cdot,\cdot]_{\triangleright},  [\cdot,\cdot,\cdot]_{\triangleleft}, [\cdot,\cdot,\cdot]_{\triangle},$ $ [\cdot,\cdot,\cdot]_{\diamond})$
and denoted by $\mathfrak{g}_\star$. The NS-3-Leibniz algebra $(\mathfrak{g},[\cdot,\cdot,\cdot]_{\triangleright},  [\cdot,\cdot,\cdot]_{\triangleleft}, [\cdot,\cdot,\cdot]_{\triangle}, [\cdot,\cdot,\cdot]_{\diamond})$
 is called a compatible NS-3-Leibniz algebra of the 3-Leibniz algebra $\mathfrak{g}_\star$.
\end{prop}

\begin{proof}
For all $ x, y, z, a, b\in \mathfrak{g}$,
by summing up the left-hand sides of the Eqs.  \eqref{7.1}- \eqref{7.6}, we can easily get  $[a,b,[x,y,z]_{\star}]_{\star}$. On the
other hand, by summing up the right-hand sides of the Eqs.  \eqref{7.1}- \eqref{7.6}, we have
$[[a,b,x]_{\star},y,z]_{\star}+[x,[a,b,y]_{\star},z]_{\star}+[x,y,[a,b,z]_{\star}]_{\star}$.
Hence the proof.
\end{proof}

Let $(\mathfrak{g},[\cdot,\cdot,\cdot]_{\triangleright},  [\cdot,\cdot,\cdot]_{\triangleleft}, [\cdot,\cdot,\cdot]_{\triangle}, [\cdot,\cdot,\cdot]_{\diamond})$  be an NS-3-Leibniz algebra.
 Define the trilinear
maps $\rho^l_{\triangleright}, \rho^m_{\triangle},$ $ \rho^r_{\triangleleft}, \Phi: \mathfrak{g}\otimes\mathfrak{g}\otimes\mathfrak{g}\rightarrow \mathfrak{g}$ by
 \begin{align*}
&\rho^l_{\triangleright}(x,y,z)=[x,y,z]_{\triangleright},~ \rho^m_{\triangle}(x,y,z)=[x,y,z]_{\triangle},~\rho^r_{\triangleleft}(x,y,z)= [x,y,z]_{\triangleleft},~\Phi(x,y,z)=[x,y,z]_{\diamond}.
\end{align*}
 With these notations, we
have the following proposition.

\begin{prop} \label{prop: representation-NS-3-Leibniz algebra}
Let $(\mathfrak{g},[\cdot,\cdot,\cdot]_{\triangleright},  [\cdot,\cdot,\cdot]_{\triangleleft}, [\cdot,\cdot,\cdot]_{\triangle}, [\cdot,\cdot,\cdot]_{\diamond})$  be an NS-3-Leibniz algebra.
 Then $(\mathfrak{g}:\rho^l_{\triangleright}, \rho^m_{\triangle}, \rho^r_{\triangleleft})$ is a representation of the
subadjacent 3-Leibniz  algebra $\mathfrak{g}_\star$ and $\Phi$ defined above is a 2-cocycle. Moreover, the identity map
$\mathrm{Id}:\mathfrak{g}\rightarrow \mathfrak{g}$ is a $\Phi$-twisted Rota-Baxter operator on the 3-Leibniz  algebra $\mathfrak{g}_\star$ with respect to $(\mathfrak{g}:\rho^l_{\triangleright}, \rho^m_{\triangle}, \rho^r_\triangleleft)$.
\end{prop}

\begin{proof}
For any $a,b,x,y,z\in    \mathfrak{g}$, by  Eq.  \eqref{7.1},  we have
\begin{align*}
\rho^l_\triangleright(a,b,\rho^l_\triangleright(x,y,z))=&[a,b,[x,y,z]_\triangleright]_\triangleright\\
=&[[a, b, x]_{\star},y,z]_{\triangleright}+ [x,  [a, b, y]_{\star},z]_{\triangleright}+ [x,y,[a, b, z]_{\triangleright}]_{\triangleright}\\
=&\rho^l_\triangleright([a,b,x]_{\star}, y, z)+ \rho^l_\triangleright(x,[a,b,y]_{\star},z)+\rho^l_\triangleright(x,y,\rho^l_\triangleright(a,b,z)).
\end{align*}
Thus, the  Eq.  \eqref{2.2}  holds. The same for Eqs. \eqref{2.3}-\eqref{2.6}. Therefore, $(\mathfrak{g}:\rho^l_{\triangleright}, \rho^m_{\triangle}, \rho^r_{\triangleleft})$ is a representation of  $\mathfrak{g}_\star$.
Moreover, Eq.  \eqref{7.6} is equivalent to $\Phi$ is a
2-cocycle  of the 3-Leibniz algebra $\mathfrak{g}_\star$ with coefficients
in the representation $(\mathfrak{g}:\rho^l_{\triangleright}, \rho^m_{\triangle}, \rho^r_\triangleleft)$. Finally, we have
\begin{align*}
&\mathrm{Id}\big(\rho^l(\mathrm{Id}x,\mathrm{Id}y,z)+\rho^m(\mathrm{Id}x,y,\mathrm{Id}z)+\rho^r(x,\mathrm{Id}y,\mathrm{Id}z)+\Phi(\mathrm{Id}x,\mathrm{Id}y,\mathrm{Id}z)\big)\\
&=[x,y,z]_\triangleleft+[x,y,z]_\triangleright+[x,y,z]_{\bigtriangleup}+[x,y,z]_\diamond\\
 &=[x,y,z]_\star=[\mathrm{Id}x,\mathrm{Id}y,\mathrm{Id}z]_\star,
\end{align*}
which shows that $\mathrm{Id}$ is a $\Phi$-twisted Rota-Baxter operator on the 3-Leibniz  algebra $\mathfrak{g}_\star$ with respect to $(\mathfrak{g}:\rho^l_{\triangleright}, \rho^m_{\triangle}, \rho^r_\triangleleft)$.
\end{proof}

\begin{theorem}\label{theorem:twisted-NS-3-Leibniz algebra}
Let $T:V\rightarrow  \mathfrak{g}$ be  a $\Phi$-twisted Rota-Baxter operator on a 3-Leibniz algebra $(\mathfrak{g}, [\cdot,\cdot,\cdot]_\mathfrak{g})$ with respect to a representation $(V;\rho^l,\rho^m,\rho^r)$.
Then there is an NS-3-Leibniz algebra
structure on $V$ with trilinear operations given by
\begin{align*}
& [u,v,w]_{\triangleleft}=\rho^r(u,Tv,Tw),~ [u,v,w]_{\triangle}=\rho^m(Tu,v,Tw),\\
& [u,v,w]_{\triangleright}=\rho^l(Tu,Tv,w),~[u,v,w]_{\diamond}=\Phi(Tu,Tv,Tw).
\end{align*}
for all $u,v,w\in   V.$
\end{theorem}
\begin{proof}
For any $u,v,w, s,t\in   V$, we have
 \begin{align*}
&[u,v,[s,t,w]_\triangleright]_{\triangleright}\\
=&\rho^l(Tu,Tv,\rho^l(Ts,Tt,w))\\
=&\rho^l([Tu,Tv, Ts]_\mathfrak{g},Tt,w)+\rho^l(Ts,[Tu,Tv,Tt]_\mathfrak{g},w)+\rho^l(Ts,Tt,\rho^l(Tu,Tv,w))\\
=&[[u, v, s]_{\star},t,w]_{\triangleright}+ [s,  [u, v, t]_{\star},w]_{\triangleright}+ [s,t,[u, v, w]_{\triangleright}]_{\triangleright}.
\end{align*}
Thus, the  Eq.  \eqref{7.1}  holds. The same for Eqs. \eqref{7.2}-\eqref{7.5}. Finally, by $\Phi$ is a 2-cocycle,  we have
 \begin{align*}
&[s,t,[u,v,w]_\diamond]_{\triangleright}+[s,t,[u,v,w]_\star]_{\diamond}+ [s,[u,v, t]_{\star},w]_{\diamond}+[[u,  v, s]_{\star},t, w]_{\diamond}+\\
&[s,  [u, v, t]_{\diamond},w]_{\triangle}+[[u,v,s]_{\diamond}, t, w]_{\triangleleft}\\
=&\rho^l(Ts,Tt,\Phi(Tu,Tv,Tw))+\Phi(Ts,Tt,[Tu,Tv,Tw]_\mathfrak{g})+\Phi(Ts,[Tu,Tv,Tt]_{\mathfrak{g}},Tw)+\\
&\Phi([Tu,Tv,Ts]_{\mathfrak{g}},Tt,Tw)+\rho^m(Ts, \Phi(Tu,Tv,Tt),Tw)+\rho^r(\Phi(Tu,Tv,Ts),Tt, Tw)\\
=&\Phi(Tu,Tv,[Ts,Tt,Tw]_{\mathfrak{g}})+\rho^l(Tu,Tv, \Phi(Ts,Tt, Tw))\\
  =& [u,v,[s, t, w]_{\star}]_{\diamond}+[u,v, [s, t, w]_{\diamond}]_{\triangleright}.
\end{align*}
This completes the proof.
\end{proof}

The subadjacent 3-Leibniz algebra of the NS-3-Leibniz algebra constructed in Proposition \ref{prop:sub-adjacent 3-Leibniz}
is given by
\begin{align*}
& [u,v,w]_{\star}=\rho^r(u,Tv,Tw)+\rho^m(Tu,v,Tw)+\rho^l(Tu,Tv,w)+\Phi(Tu,Tv,Tw).
\end{align*} for all $u,v,w\in V.$
This 3-Leibniz algebra structure on $V$ coincides with the one given in Eq.  \eqref{3.2}.

By combining Example \ref{exam:Nijenhuis} and Theorem \ref{theorem:twisted-NS-3-Leibniz algebra}, we can get that a Nijenhuis operator on a 3-Leibniz algebra induces an NS-3-Leibniz algebra.

\begin{exam}
Let $N:\mathfrak{g}\rightarrow\mathfrak{g}$ be a Nijenhuis operator on a 3-Leibniz algebra $(\mathfrak{g}, [\cdot,\cdot,\cdot]_\mathfrak{g})$.
Then the trilinear operations
\begin{align*}
& [x,y,z]_{\triangleleft}=[x,Ny,Nz]_\mathfrak{g},~ [x,y,z]_{\triangle}=[Nx,y,Nz]_\mathfrak{g},~ [x,y,z]_{\triangleright}=[Nx,Ny,z]_\mathfrak{g},\\
&[x,y,z]_{\diamond}=-N\big([x,y,Nz]_\mathfrak{g}+[x,Ny,z]_\mathfrak{g}+[Nx,y,z]_\mathfrak{g}\big)+N^2[x,y,z]_\mathfrak{g}
\end{align*}
defines an NS-3-Leibniz algebra structure on $\mathfrak{g}$.
\end{exam}

Combined with Example \ref{exam:Reynolds} and Theorem \ref{theorem:twisted-NS-3-Leibniz algebra},   a Reynolds operator on a 3-Leibniz algebra also induces an NS-3-Leibniz algebra.

\begin{exam}
Let $T:\mathfrak{g}\rightarrow\mathfrak{g}$ be a  Reynolds operator on a 3-Leibniz algebra $(\mathfrak{g}, [\cdot,\cdot,\cdot]_\mathfrak{g})$.
Then the trilinear operations
\begin{align*}
& [x,y,z]_{\triangleleft}=[x,Ty,Tz]_\mathfrak{g},~ [x,y,z]_{\triangle}=[Tx,y,Tz]_\mathfrak{g},~ [x,y,z]_{\triangleright}=[Tx,Ty,z]_\mathfrak{g},~[x,y,z]_{\diamond}=- [Tx,Ty,Tz]_\mathfrak{g}
\end{align*}
defines an NS-3-Leibniz algebra structure on $\mathfrak{g}$.
\end{exam}

NS-3-Leibniz algebras  also arise from weighted Rota-Baxter operators on 3-Leibniz algebras. Let $(\mathfrak{g}, [\cdot,\cdot,\cdot]_\mathfrak{g})$
be a 3-Leibniz algebra. A linear map $B:\mathfrak{g}\rightarrow\mathfrak{g}$ is said to be a Rota-Baxter operator of weight $\lambda$ on the
3-Leibniz algebra if $B$ satisfies
\begin{align*}
[Bx,By,Bz]_\mathfrak{g}= &B\big([x,By,Bz]_\mathfrak{g}+[Bx,y,Bz]_\mathfrak{g}+[Bx,By,z]_\mathfrak{g}+\\
&\lambda[x,y,Bz]_\mathfrak{g}+\lambda[Bx,y,z]_\mathfrak{g}+\lambda[x,By,z]_\mathfrak{g}+\lambda^2[x,y,z]_\mathfrak{g}\big)
\end{align*}
In the following, we show that a Rota-Baxter operators of weight $\lambda$ induces an NS-3-Leibniz algebra.

\begin{prop}
Let $B:\mathfrak{g}\rightarrow\mathfrak{g}$  be a Rota-Baxter operator of weight  $\lambda$  on the 3-Leibniz algebra $(\mathfrak{g}, [\cdot,\cdot,\cdot]_\mathfrak{g})$.
Then there is an NS-3-Leibniz algebra structure on the vector space $\mathfrak{g}$ with trilinear operations
\begin{align*}
& [x,y,z]_{\triangleleft}=[x,By,Bz]_\mathfrak{g},~ [x,y,z]_{\triangle}=[Bx,y,Bz]_\mathfrak{g},~ [x,y,z]_{\triangleright}=[Bx,By,z]_\mathfrak{g},\\
&[x,y,z]_{\diamond}=\lambda[x,y,Bz]_\mathfrak{g}+\lambda[Bx,y,z]_\mathfrak{g}+\lambda[x,By,z]_\mathfrak{g}+\lambda^2[x,y,z]_\mathfrak{g}.
\end{align*}
for all $x,y,z\in \mathfrak{g}.$
\end{prop}

\begin{proof}
For any $x,y,z, a,b\in \mathfrak{g}$, we have
 \begin{align*}
   &[[a, b, x]_{\star},y,z]_{\triangleright}+ [x,  [a, b, y]_{\star},z]_{\triangleright}+ [x,y,[a, b, z]_{\triangleright}]_{\triangleright}\\
 =&[[Ba,Bb,Bx]_\mathfrak{g},By, z]_{\mathfrak{g}}+ [Bx,  [Ba,Bb,By]_\mathfrak{g},z]_{_\mathfrak{g}}+ [Bx,By,[Ba,Bb,z]_\mathfrak{g}]_{\mathfrak{g}}\\
 =&[Ba,Bb,[Bx,By, z]_\mathfrak{g}]_{\mathfrak{g}}\\
 =&[a,b,[x,y,z]_\triangleright]_{\triangleright}
\end{align*}
 and
\begin{align*}
 &[[a, b, x]_{\triangleleft},y,z]_{\triangleleft}+ [x,  [a, b, y]_{\triangleleft},z]_{\triangle}+ [x,y,[a, b, z]_{\triangleleft}]_{\triangleright}\\
 =&[[a,Bb,Bx]_\mathfrak{g},By,B z]_{\mathfrak{g}}+ [Bx,  [a,Bb,By]_\mathfrak{g},Bz]_{_\mathfrak{g}}+ [Bx,By,[a,Bb,Bz]_\mathfrak{g}]_{\mathfrak{g}}\\
 =&[a,Bb,[Bx,By,B z]_\mathfrak{g}]_{\mathfrak{g}}=[a,Bb,B[x,y,z]_\star]_{\mathfrak{g}}\\
 =&[a,b,[x,y,z]_\star]_{\triangleleft}.
\end{align*}
Also,
\begin{align*}
 &[[a, b, x]_{\triangle},y,z]_{\triangleleft}+ [x,  [a, b, y]_{\triangle},z]_{\triangle}+ [x,y,[a, b, z]_{\triangle}]_{\triangleright}\\
 =&[[Ba,b,Bx]_\mathfrak{g},By,B z]_{\mathfrak{g}}+ [Bx,  [Ba,b,By]_\mathfrak{g},Bz]_{_\mathfrak{g}}+ [Bx,By,[Ba,b,Bz]_\mathfrak{g}]_{\mathfrak{g}}\\
 =&[Ba,b,[Bx,By,B z]_\mathfrak{g}]_{\mathfrak{g}}=[Ba,b,B[x,y,z]_\star]_{\mathfrak{g}}\\
 =&[a,b,[x,y,z]_\star]_{\triangle}.
\end{align*}
Similarly, we have
\begin{align*}
 &[[a, b, x]_{\triangleright},y,z]_{\triangleleft}+ [x,  [a, b, y]_{\star},z]_{\triangleleft}+ [x,y,[a, b, z]_{\star}]_{\triangleleft}\\
 =&[[Ba,Bb,x]_\mathfrak{g},By,B z]_{\mathfrak{g}}+ [x,  [Ba,Bb,By]_\mathfrak{g},Bz]_{_\mathfrak{g}}+ [x,By,[Ba,Bb,Bz]_\mathfrak{g}]_{\mathfrak{g}}\\
 =&[Ba,Bb,[x,By,B z]_\mathfrak{g}]_{\mathfrak{g}}\\
 =&[a,b,[x,y,z]_\triangleleft]_{\triangleright},\\
  &[[a, b, x]_{\star},y,z]_{\triangle}+ [x,  [a, b, y]_{\triangleright},z]_{\triangle}+ [x,y,[a, b, z]_{\star}]_{\triangle}\\
 =&[[Ba,Bb,Bx]_\mathfrak{g},y,B z]_{\mathfrak{g}}+ [Bx,  [Ba,Bb,y]_\mathfrak{g},Bz]_{_\mathfrak{g}}+ [Bx,y,[Ba,Bb,Bz]_\mathfrak{g}]_{\mathfrak{g}}\\
 =&[Ba,Bb,[Bx,y,B z]_\mathfrak{g}]_{\mathfrak{g}}\\
 =&[a,b,[x,y,z]_\diamond]_{\triangleright}.
 \end{align*}
Finally, we have
 \begin{align*}
   &[a,b,[x, y, z]_{\star}]_{\diamond}+[a,b,[x, y, z]_{\diamond}]_{\triangleright}\\
=&\lambda[a,b,[B x, By, Bz]_{ \mathfrak{g}}]_\mathfrak{g}+\lambda[Ba,b,[x, y, z]_{\star}]_\mathfrak{g}+\lambda[a,Bb,[x, y, z]_{\star}]_\mathfrak{g}+\\
&\lambda^2[a,b,[x, y, z]_{\star}]_\mathfrak{g}+[Ba,Bb,[x, y, z]_{\diamond}]_{\mathfrak{g}}\\
=&\lambda[a,b,[B x, By, Bz]_{ \mathfrak{g}}]_\mathfrak{g}+\lambda[Ba,b,[x,By,Bz]_\mathfrak{g}]_\mathfrak{g}+\lambda[Ba,b,[Bx,y,Bz]_\mathfrak{g}]_\mathfrak{g}+\lambda[Ba,b,[Bx,By,z]_\mathfrak{g}]_\mathfrak{g}+\\
&\lambda^2[Ba,b,[x,y,Bz]_\mathfrak{g}]_\mathfrak{g}+\lambda^2[Ba,b,[Bx,y,z]_\mathfrak{g}]_\mathfrak{g}+\lambda^2[Ba,b,[x,By,z]_\mathfrak{g}]_\mathfrak{g}+\lambda^3[Ba,b,[x,y,z]_\mathfrak{g}]_\mathfrak{g}+\\
&\lambda[a,Bb,[x,By,Bz]_\mathfrak{g}]_\mathfrak{g}+\lambda[a,Bb,[Bx,y,Bz]_\mathfrak{g}]_\mathfrak{g}+\lambda[a,Bb,[Bx,By,z]_\mathfrak{g}]_\mathfrak{g}+\lambda^2[a,Bb,[x,y,Bz]_\mathfrak{g}]_\mathfrak{g}+\\
&\lambda^2[a,Bb,[Bx,y,z]_\mathfrak{g}]_\mathfrak{g}+\lambda^2[a,Bb,[x,By,z]_\mathfrak{g}]_\mathfrak{g}+\lambda^3[a,Bb,[x,y,z]_\mathfrak{g}]_\mathfrak{g}+\\
&\lambda^2 [a, b,[x,By,Bz]_\mathfrak{g}]_\mathfrak{g}+\lambda^2[a, b,[Bx,y,Bz]_\mathfrak{g}]_\mathfrak{g}+\lambda^2[a, b,[Bx,By,z]_\mathfrak{g}]_\mathfrak{g}+\lambda^3[a, b,[x,y,Bz]_\mathfrak{g}]_\mathfrak{g}+\\
&\lambda^3 [a, b,[Bx,y,z]_\mathfrak{g}]_\mathfrak{g}+\lambda^3[a, b,[x,By,z]_\mathfrak{g}]_\mathfrak{g}+\lambda^4[a, b,[x,y,z]_\mathfrak{g}]_\mathfrak{g}+\\
&\lambda[Ba,Bb,[x,y,Bz]_\mathfrak{g}]_{\mathfrak{g}}+\lambda[Ba,Bb,[Bx,y,z]_\mathfrak{g}]_{\mathfrak{g}}+\lambda[Ba,Bb,[x,By,z]_\mathfrak{g}]_{\mathfrak{g}}+\lambda^2[Ba,Bb,[x,y,z]_\mathfrak{g}]_{\mathfrak{g}}\\
=&[x,y,[a,b,z]_\diamond]_{\triangleright}+[x,y,[a,b,z]_\star]_{\diamond}+[x,  [a,b,y]_{\diamond},z]_{\triangle}+ [x,[a,b, y]_{\star},z]_{\diamond}+\\
&[[a,b,x]_{\star},y, z]_{\diamond}+[[a,b,x]_{\diamond},y, z]_{\triangleleft}.
\end{align*}
The proof is finished.
\end{proof}

Below we give a necessary and sufficient condition for the existence of a compatible NS-3-Leibniz
algebra structure on a 3-Leibniz algebra.

\begin{prop}
Let  $(\mathfrak{g}, [\cdot,\cdot,\cdot]_\mathfrak{g})$  be a 3-Leibniz algebra.
 Then there is a compatible NS-3-Leibniz algebra
structure on $\mathfrak{g}$ if and only if there exists an invertible $\Phi$-twisted Rota-Baxter operator
$T:V\rightarrow \mathfrak{g}$ on $\mathfrak{g}$
with respect to a representation $(V;\rho^l,\rho^m,\rho^r)$ and a 2-cocycle $\Phi$. Furthermore, the compatible NS-3-Leibniz
algebra structure on $\mathfrak{g}$ is given by
\begin{align*}
& [x,y,z]_{\triangleleft}=T\rho^r(T^{-1}x,y,z),~ [x,y,z]_{\triangle}=T\rho^m(x,T^{-1}y,z),\\
& [x,y,z]_{\triangleright}=T\rho^l(x,y,T^{-1}z),~[x,y,z]_{\diamond}=T\Phi(x,y,z).
\end{align*}
\end{prop}

\begin{proof}
Let $T:V\rightarrow \mathfrak{g}$ be an invertible $\Phi$-twisted Rota-Baxter operator on $\mathfrak{g}$
with respect to a representation $(V;\rho^l,\rho^m,\rho^r)$ and a 2-cocycle $\Phi$.
By  Theorem \ref{theorem:twisted-NS-3-Leibniz algebra}, there is a NS-3-Leibniz algebra structure on $V$ given by
\begin{align*}
& [u,v,w]_{\triangleleft}=\rho^r(u,Tv,Tw),~ [u,v,w]_{\triangle}=\rho^m(Tu,v,Tw),\\
& [u,v,w]_{\triangleright}=\rho^l(Tu,Tv,w),~[u,v,w]_{\diamond}=\Phi(Tu,Tv,Tw).
\end{align*}
Since $T$ is an invertible map, the trilinear operations
\begin{align*}
& [x,y,z]_{\triangleleft}=T[T^{-1}x,T^{-1}y,T^{-1}z]_{\triangleleft}=T\rho^r(T^{-1}x,y,z),\\
& [x,y,z]_{\triangle}=T[T^{-1}x,T^{-1}y,T^{-1}z]_{\triangle}=T\rho^m(x,T^{-1}y,z),\\
& [x,y,z]_{\triangleright}=T[T^{-1}x,T^{-1}y,T^{-1}z]_{\triangleright}=T\rho^l(x,y,T^{-1}z),\\
& [x,y,z]_{\diamond}=T[T^{-1}x,T^{-1}y,T^{-1}z]_{\diamond}=T\Phi(x,y,z).
\end{align*}
defines a NS-3-Leibniz algebra on $\mathfrak{g}$.
Moreover, we have
 \begin{align*}
&[x,y,z]_\triangleleft+[x,y,z]_\triangleright+[x,y,z]_{\bigtriangleup}+[x,y,z]_\diamond\\
=&T\rho^r(T^{-1}x,y,z)+T\rho^l(x,y,T^{-1}z)+T\rho^m(x,T^{-1}y,z)+T\Phi(x,y,z)\\
=&T\big(\rho^r(T^{-1}x,T\circ T^{-1}y,T\circ T^{-1}z)+\rho^l(T\circ T^{-1}x,T\circ T^{-1}y,T^{-1}z)+\rho^m(T\circ T^{-1}x,T^{-1}y,T\circ T^{-1}z)+\\
&\Phi(T\circ T^{-1}x,T\circ T^{-1}y,T\circ T^{-1}z)\big)\\
=&[T\circ T^{-1} x,T\circ T^{-1}  y,  T \circ T^{-1} z]_\mathfrak{g}=[x,y,z]_\mathfrak{g}.
\end{align*}

Conversely, let $(\mathfrak{g},[\cdot,\cdot,\cdot]_{\triangleright},[\cdot,\cdot,\cdot]_{\triangleleft}, [\cdot,\cdot,\cdot]_{\triangle}, [\cdot,\cdot,\cdot]_{\diamond})$
 be a compatible NS-3-Leibniz algebra structure on $\mathfrak{g}$. By Proposition \ref{prop: representation-NS-3-Leibniz algebra},
$(\mathfrak{g}:\rho^l_{\triangleright}, \rho^m_{\triangle}, \rho^r_\triangleleft)$ is a representation of the Leibniz algebra $(\mathfrak{g},[\cdot,\cdot,\cdot]_\mathfrak{g})$,
and the identity map
$\mathrm{Id}:\mathfrak{g}\rightarrow \mathfrak{g}$ is a $\Phi$-twisted Rota-Baxter operator on the 3-Leibniz  algebra $(\mathfrak{g},[\cdot,\cdot,\cdot]_\mathfrak{g})$ with respect to $(\mathfrak{g}:\rho^l_{\triangleright}, \rho^m_{\triangle}, \rho^r_\triangleleft)$.
\end{proof}

\noindent
{{\bf Acknowledgments.} The work is supported by  the Science and Technology Program of Guizhou Province (Grant No.\, QKHJC QN[2025]362).}

\end{document}